\documentclass[reqno,11pt]{amsart}
\usepackage{a4wide,xcolor,eucal,enumerate,mathrsfs}
\usepackage[normalem]{ulem}
\usepackage{pdfsync}
\usepackage[hidelinks]{hyperref}
\usepackage{amsmath,amssymb,epsfig,amsthm}
\usepackage[latin1]{inputenc}

\numberwithin{equation}{section}

\newtheorem{theorem}{Theorem}[section]

\newtheorem{lemma}[theorem]{Lemma}
\newtheorem{proposition}[theorem]{Proposition}
\newtheorem{definition}[theorem]{Definition}

\newtheorem{example}[theorem]{Example}
\newtheorem{remark}[theorem]{Remark}

\newcommand{\N}{\mathbb{N}}
\newcommand{\Q}{\mathbb{Q}}
\newcommand{\R}{\mathbb{R}}

\newcommand{\V}{\mathbb{V}}

\newcommand{\cA}{{\ensuremath{\mathcal A}}}

\newcommand{\cE}{{\ensuremath{\mathcal E}}}

\newcommand{\cL}{{\ensuremath{\mathcal L}}}

\renewcommand{\aa}{{\mbox{\boldmath$a$}}}
\newcommand{\bb}{{\mbox{\boldmath$b$}}}
\newcommand{\cc}{{\mbox{\boldmath$c$}}}

\newcommand{\mm}{{\mbox{\boldmath$m$}}}

\newcommand{\XX}{{\mbox{\boldmath$X$}}}
\newcommand{\YY}{{\mbox{\boldmath$Y$}}}
\newcommand{\sXX}{{\mbox{\scriptsize\boldmath$X$}}}
\newcommand{\sYY}{{\mbox{\scriptsize\boldmath$Y$}}}

\newcommand{\sbb}{{\mbox{\scriptsize\boldmath$b$}}}

\newcommand{\eeta}{{\mbox{\boldmath$\eta$}}}
\newcommand{\llambda}{{\mbox{\boldmath$\lambda$}}}

\newcommand{\ssigma}{{\mbox{\boldmath$\sigma$}}}

\newcommand{\supp}{\mathop{\rm supp}\nolimits}   
\newcommand{\restr}[1]{\lower3pt\hbox{$|_{#1}$}}

\newcommand{\Leb}[1]{{\mathscr L}^{#1}}      
\newcommand{\BorelSets}[1]{{\mathscr B}(#1)}
\newcommand{\eps}{\varepsilon}  
\newcommand{\nchi}{{\raise.3ex\hbox{$\chi$}}}

\def\qed{\ifmmode 
  \else \leavevmode\unskip\penalty9999 \hbox{}\nobreak\hfill
  \fi               
    \qquad           \hbox{\hskip.5em $\square$
                \hskip.1em}}

\renewcommand{\mm}{\mathfrak m} 
\newcommand{\Gq}[1]{\Gamma\!\left(#1\right)}          

\newcommand{\Probabilities}[1]{\mathscr P(#1)}          

\newcommand{\res}{\mathop{\hbox{\vrule height 7pt width .5pt depth 0pt
\vrule height .5pt width 6pt depth 0pt}}\nolimits}

\newcommand{\sqGq}[1]{\sqrt{ \Gq {#1}}}
\newcommand{\bra}[1]{\left( #1 \right)}
\newcommand{\sqa}[1]{\left[ #1 \right]}
\newcommand{\cur}[1]{\left\{ #1 \right\}}

\newcommand{\abs}[1]{\left| #1 \right|}
\newcommand{\nor}[1]{\left\| #1 \right\|}
\newcommand{\Lbm}[1]{L^{#1}(\mm)}

\newcommand{\vphi}{\varphi}

\renewcommand{\div}{\mathop{\rm div}\nolimits}        

\newcommand{\dbneg}{\div\bb^-}
\newcommand{\veps}{\varepsilon}

\newcommand{\Commutator}{{\mathscr C}}
\newcommand{\Algebra}{{\mathscr A}}
\newcommand{\scrL}{\mathscr{L}}
\renewcommand{\cL}{\mathsf{L}}

\title[Lecture notes on the DiPerna-Lions theory in abstract measure spaces]{Lecture notes on the DiPerna-Lions \\ theory in abstract measure spaces}

\begin{document}

\author{Luigi Ambrosio}
   \address{Scuola Normale Superiore, Pisa} 
   \email{luigi.ambrosio@sns.it}
   \author{Dario Trevisan}
   \address{Universit\`a di Pisa} 
   \email{dario.trevisan@unipi.it}

\maketitle

\tableofcontents

\section{Introduction}

In these notes, which correspond closely to the lectures given by the first author in Toulouse  (with the exception of the
last part, dealing with more recent developments), we describe the extension of the theory of ODE well-posedness to abstract
spaces, recently developed in \cite{Ambrosio-Trevisan}. Even though the concept of ODE requires only the differential structure, 
the mathematical tools needed for uniqueness and the applications require a
``metric'' structure. Hence, the results described in these notes could be read either at the level of Dirichlet spaces or at the
level of metric measure spaces. We have chosen, as in \cite{Ambrosio-Trevisan}, to adopt the language of Dirichlet spaces, since
several tools rely on the well-established theory of diffusion semigroups and, in particular, on $\Gamma$-calculus techniques.

In order to simplify as much as possible the exposition, we have made some simplifying assumptions, considering for instance
only probability measures as reference measures and choosing model cases. Indeed, many different assumptions on the integrability
of the vector field, on its divergence, on its symmetric derivative and on the solution can be made, and we refer to \cite{Ambrosio-Trevisan}
for more general statements and detailed computations.

In order to describe the origins of the theory, let us consider the system of ODE's
\begin{equation}\label{eq:ODE-system}
\begin{cases} x'(t)=\bb(t,x(t)) &\text{$t\in (0,T)$}\\
x(0)=x,
\end{cases}
\end{equation}
where $\bb(t,x)=\bb_t(x)$ is a (possibly) time-dependent family of vector fields in $\R^n$. A natural question is the meaning that one has
to give to the ODE in the case when $\bb_t$ are not smooth, or even are defined up to Lebesgue negligible sets. This problem arises in a natural
way in fluid mechanics (where, in many cases, the velocity belongs to a Sobolev class) and in the theory of conservation laws (where
the velocity can be $BV$, or even worse), see the lecture notes \cite{cetraro}, \cite{bologna} for a more detailed account of these
developments. 

The classical well-posedness theory of \eqref{eq:ODE-system} requires the regularity
$$
{\rm Lip}(\bb_t)\in L^1(0,T)
$$
and it ensures Lipschitz regularity w.r.t.\ $x$ of the flow map $\XX(t,x)$, as well as stability w.r.t.\ approximations of
the vector field. For vector fields with a special structure, for instance satisfying a one-sided Lipschitz condition
$$
\langle \bb_t(x)-\bb_t(y),x-y\rangle\leq C(t)|x-y|^2\qquad C\in L^1(0,T)
$$
or autonomous and gradient vector fields $\bb=-\nabla V$ with $V$ convex, still pointwise uniqueness and stability hold.
More generally, under a local Lipschitz
condition one obtains a unique maximal flow, with a lower semicontinuous maximal existence time $T_\sbb:\R^n\to (0,T]$.

In \cite{DiPerna-Lions89}, in a seminal paper, DiPerna and Lions first introduced an appropriate notion of ``almost everywhere well posedness''
suitable for nonsmooth (and, in particular, Sobolev) vector fields. In \cite{Ambrosio03}, the first author revisited and improved the axiomatization 
of \cite{DiPerna-Lions89}, emphasizing in particular the connection with probability and with the theory of stochastic processes. In order to illustrate this
point of view, let us start from the classical nonuniqueness example on the real line:
$$
x'=\sqrt{|x|},\qquad x_0=-c^2,\,\,c\geq 0.
$$
In this case one has $x(t)=-(t/2-c)^2$ for $0\leq t\leq 2c$, then the solution can stay at the origin for some time $2T(c)$ (possibly infinite, or null)
and then continue as $x(t)=(t/2-T(c)-c)^2$ for $t\geq 2c+2T(c)$, if $T(c)<\infty$. Assume that we choose (measurably) the time $T(c)$, so that we have
a flow map 
$$
\XX(t,x):[0,\infty)\times (-\infty,0]\to\R.
$$
Simple constructions show that this flow map is very sensitive to the approximation of $\bb$, namely different (Lipschitz) approximations
of $\bb$ produce in the limit different flow maps $\XX$. So, what is special about the ``natural'' solution with $T(c)=0$? It is simple
to check that
$$
T(c)=0\quad\text{for a.e.\ $c\leq 0$}\qquad\Longleftrightarrow\qquad
\XX(t,\cdot)_\#\Leb{1}\ll\Leb{1}\quad\text{for a.e.\ $t\in (0,\infty)$},
$$
where we use the standard notation for the push-forward operator between measures. Indeed, $\XX(t,\cdot)_\#\Leb{1}(0)=0$
for a.e.\ $t\in (0,\infty)$ implies that the set
$$
\{(t,c):\ \XX(t,-c^2)=0\}
$$
is negligible w.r.t.\ $\Leb{1}\times\Leb{1}$, hence Fubini's theorem gives that for $\Leb{1}$-a.e.\ $c\in (-\infty,0]$ one has $T(c)=0$.

If we make the absolute continuity condition quantitative, this leads to the concept of Regular Lagrangian Flow.

\begin{definition}[Regular Lagrangian Flow] We say that $\XX:[0,T]\times\R^n\to\R^n$ is a Regular Lagrangian flow if:
\begin{itemize}
\item[(a)] for $\Leb{n}$-a.e.\ $x\in\R^n$ the function $t\mapsto\XX(t,x)$ is an absolutely continuous solution to the
Cauchy problem \eqref{eq:ODE-system};
\item[(b)] there exists $C\in [0,\infty)$ satisfying $\XX(t,\cdot)_\#\Leb{n}\leq C\Leb{n}$ for all $t\in [0,T]$.
\end{itemize}
\end{definition}

Again, a simple application of Fubini's theorem
shows that for any Lebesgue negligible Borel set $N\subset (0,T)\times\R^n$, the set
$$
\bigl\{t \in (0,T):\ (t,\XX(t,x))\in N\bigr\}
$$
is $\Leb{1}$-negligible for $\Leb{n}$-a.e.\ $x\in\R^n$. Hence, under condition (b), the seemingly unstable condition (a)
is invariant under modifications of $\bb$ in $\Leb{1+n}$-negligible sets.

Notice also that the second condition is a global one: we are selecting a family of ODE solutions in such a way that for all
times the trajectories do not concentrate too much, relative to $\Leb{n}$.
Because of this global aspect, it is a bit
misleading to think that uniqueness of RLF implies ``existence of a unique solution starting from a.e.\ initial point''.
Nevertheless, a basic output of the DiPerna-Lions theory is the following existence and uniqueness result.

\begin{theorem} \label{thm:maindp} Let $T\in (0,\infty)$. Assume that 
$$
\bb\in L^1\bigl((0,T); W^{1,1}_{\rm loc}(\R^n;\R^n)\bigr),\qquad
\div\bb_t\in L^1\bigl((0,T);L^\infty(\R^n)\bigr),
$$
and that
\begin{equation}\label{eq:growth}
\frac{|\bb|}{1+|x|}\in L^1\bigl((0,T);L^1(\R^n)\bigr)+L^1\bigl((0,T);L^\infty(\R^n)\bigr).
\end{equation}
Then, there exists a unique RLF associated to $\bb$ (more precisely, if $\XX$ and $\YY$ are RLF's, then
$\XX(\cdot,x)=\YY(\cdot,x)$ for $\Leb{n}$-a.e.\ $x\in\R^n$).
\end{theorem}

The goal of these lectures is to provide a reasonable extension of the previous theorem to abstract spaces.
Of course we have to specify what we mean by:

\smallskip\noindent
$\bullet$ Vector field and solution to the ODE

\smallskip\noindent
$\bullet$ Sobolev space of vector fields and divergence

\smallskip\noindent
in our abstract setup and we achieve this goal using the formalism of $\Gamma$-calculus. 
After doing this, we shall provide existence and uniqueness results, both at the Lagrangian and the
Eulerian level, very similar to those of the Euclidean theory.

In the last section of the notes we describe, without proofs, further developments contained in the second author's PhD thesis
\cite{Trevisan-14}. Still following the point of view of $\Gamma$-calculus one can define a diffusion operator $\cL$ by the
requiring the validity of the formula
\[ \frac 1 2 \sqa{ \cL(f g ) - f \cL(g) - g \cL(f)} = \aa(\nabla f, \nabla g),\]
where $\aa$ is a bilinear operator satisfying the Leibniz rule. Building on this, the ``deterministic'' results presented in the
main part of these lectures can be viewed as particular cases (or limiting cases) when the diffusion coefficient $\aa$ is $0$, i.e.\ the
operator $\cL$ contains only a transport term.

\section{Reminders on the Cauchy-Lipschitz theory}

Via the theory of characteristics, the flow map $\XX$ provides also the unique solutions to the continuity
equation and the transport equation. The continuity equation takes the conservative form
\begin{equation}\label{eq:CE}
\frac{d}{dt}u_t+\div(\bb_tu_t)=0
\end{equation}
and the solution is provided, in measure-theoretic terms, by 
\begin{equation}\label{eq:1}
u_t\Leb{n}=\XX(t,\cdot)_\#(\bar u\Leb{n})
\end{equation}
see also \eqref{eq:2} for an explicit formula, via change of variables. The transport equation
\begin{equation}\label{eq:TE}
\frac{d}{dt}w_t+\bb_t\cdot\nabla w_t=c_tw_t
\end{equation}
can be solved by integrating the linear ODE
\begin{equation*} 
\frac{d}{dt}w_t(\XX(t,x))=c_t(\XX(t,x))w_t(\XX(t,x)).
\end{equation*}

In the abstract setup, the continuity equation will be more useful than the transport equation and so
we will focus mostly on it, even though the techniques developed in \cite{Ambrosio-Trevisan} can be
used to obtain well-posedness results also for general transport equations. Notice also the continuity
equation, because of its conservative form, has a weak formulation, for instance writing
$$
\frac{d}{dt}\int_{\R^n}f u_t\,dx=\int_{\R^n}\bb_t\cdot\nabla f\,  u_t\,dx
\qquad\text{in $\mathcal D'(0,T)$}
$$
for all $f\in C^1_c(\R^n)$. For the weak
formulation of \eqref{eq:TE} we need that $\div\bb_t$ is a function, writing it in the form
$$
\frac{d}{dt}w_t+\div(\bb_t w_t)=(c_t+\div\bb_t) w_t,
$$
i.e., in the weak formulation,
$$
\frac{d}{dt}\int_{\R^n}f w_t\,dx=\int_{\R^n}\bigl(\bb_t\cdot\nabla f\,  w_t+(c_t+\div\bb_t)w_tf\bigr)\,dx
\qquad\text{in $\mathcal D'(0,T)$}
$$
for all $f\in C^1_c(\R^n)$. This shows the equivalence of the two equations when $c_t=-\div\bb_t$.

\section{Nonsmooth vector fields in Euclidean spaces}

In this section we illustrate the basic ingredients of the proof of Theorem~\ref{thm:maindp}. First of all, there is a general transfer
mechanism for well posedness which allows to deduce uniqueness of RLF's from uniqueness of solutions to the continuity equation
\eqref{eq:CE}.
This transfer mechanism will be described in the abstract setup, so let us focus on uniqueness, treating also existence just
for the sake of completeness.

The proof of existence requires no regularity on $\bb$, besides the growth bounds on $\bb$ and on its divergence: indeed, if 
we mollify the vector field w.r.t.\ the space variable, we have that solutions $u^\eps$ to 
$$
\frac{d}{dt}v_t+\div (\bb^\eps_tv_t)=0
$$
with the initial condition $\bar u$ can be derived out of \eqref{eq:1} with the change of variables. It turns out that
they are provided by the explicit formula
\begin{equation}\label{eq:2}
u^\eps(t,x)\circ\XX^\eps(t,x)=\frac{\bar u(x)}{J\XX^\eps(t,x)},
\end{equation}
where $\XX^\eps$ is the flow associated to $\bb^\eps$ and $J\XX^\eps(t,x)={\rm det}\bigl(\nabla_x\XX^\eps(t,x))$.
On the other hand, the $L^1(L^\infty)$ lower bound on the divergence on $\bb_t$ is retained by $\bb^\eps$, hence the
Liouville ODE
$$
\frac{d}{dt}J\XX^\eps(t,x)=\div\bb^\eps_t(\XX^\eps(t,x))J\XX^\eps(t,x),\qquad J\XX^\eps(0,x)=1
$$
ensures that $(u^\eps)$ is uniformly bounded in $L^\infty(L^\infty)$. By the change of variables formula, $(u^\eps)$
is uniformly bounded in $L^\infty(L^1)$ as well, so taking $w^*$-limits as $\eps\searrow 0$ provides the following
result:

\begin{proposition} Under the assumptions of Theorem~\ref{thm:maindp}, the continuity equation \eqref{eq:CE} has a solution
in $L^\infty(L^1\cap L^\infty)$ for any initial datum $\bar u\in L^1\cap L^\infty$.
\end{proposition}

Notice that for the validity of this argument only $L^\infty$ lower bounds on $\div\bb_t$ are needed; this improvement has been
noticed in \cite{Ambrosio03} and it will persist also in the theory in general spaces. Also, the Sobolev regularity did not play any
role in the proof of the proposition. However, this regularity will play a decisive role in the proof of uniqueness: without this assumption (see
for instance \cite{depauw}), uniqueness may fail even for divergence-free vector fields. It turns out that, if $\bb_t$ is sufficiently
regular, weak solutions to the continuity (or transport) equation enjoy a crucial regularity property.

\begin{definition} [Renormalized solution] We say that a weak solution $u\in L^\infty(L^1\cap L^\infty)$ to \eqref{eq:CE} is renormalized if
\begin{equation}\label{eq:4}
\frac{d}{dt}\beta(u_t)+\bb_t\cdot\nabla \beta(u_t)=-u_t\beta'(u_t)\div\bb_t,
\end{equation}
still in the weak sense, for all $\beta\in C^1(\R)$. 
\end{definition}

Rewriting \eqref{eq:4} in conservative form, we have
$$
\frac{d}{dt}\beta(u_t)+\div(\bb_t\beta(u_t))=(\beta(u_t)-u_t\beta'(u_t))\div\bb_t,
$$
so that the weak formulation of \eqref{eq:4} is
$$
\frac{d}{dt}\int_{\R^n}\phi\beta(u_t)\,dx=\int_{\R^n}\beta(u_t) \bb_t \cdot\nabla\phi\,dx+\int_{\R^n}\phi\bigl(\beta(u_t)-u_t\beta'(u_t)\bigr)\div\bb_t\,dx
$$
for all $\phi\in C^1_c(\R^n)$.

Then, well posedness of \eqref{eq:CE} is a consequence of the following result.

\begin{theorem} Assume that $\div\bb_t\in L^\infty(L^\infty)$ and that
all weak solutions $u\in L^\infty(L^1\cap L^\infty)$ to \eqref{eq:TE} with $c\in L^\infty(L^\infty)$
are renormalized. Under the growth assumption \eqref{eq:growth} on $\bb$, the PDE \eqref{eq:CE} is
well posed in $L^\infty(L^1\cap L^\infty)$.
\end{theorem}

Let us give a sketchy proof of this result. We already know about existence, so we have to prove uniqueness.
By the linearity of the equation, we need to show that $u\equiv 0$ if the initial datum $\bar u$ is $0$. First, we
extend $\bb$ to negative times, setting $\bb_t=0$ for all $t<0$, and we extend $u$ accordingly, setting
$u_t=0$ for all $t<0$. Then, it is easily seen that $u$ is a weak solution in $(-\infty,T)\times\R^n$, so we
can apply the renormalization property with $\beta_\epsilon(z)=\sqrt{z^2+\epsilon^2}-\epsilon$ to get
$$
\frac{d}{dt}\int_{\R^n}\beta_\eps(u_t)\phi\,dx=\int_{\R^n}\beta_\eps(u_t) \bb_t \cdot \nabla\phi dx+
\int_{\R^n}\phi\bigl(\beta_\eps(u_t)-u_t\beta_\eps'(u_t)\bigr)\div\bb_t\,dx
$$
in the sense of distributions in $(-\infty,T)$, for all $\phi\in C^1_c(\R^n)$. Since
$\bigl|\beta_\eps(u_t)-u_t\beta_\eps'(u_t)\bigr|\leq\eps$ one gets
$$
\frac{d}{dt}\int_{\R^n}\beta_\eps(u_t)\phi\,dx=\int_{\R^n}\beta_\eps(u_t) \bb_t \cdot \nabla\phi dx+
\eps\int_{\R^n}\phi|\div\bb_t|\,dx
$$
We let first $\eps\searrow 0$ and then use a family of cut-off functions $\phi=\phi_R\uparrow 1$
(and in this step the global growth conditions on $\bb$ and $u$ come into play), to prove that
$$
\frac{d}{dt}\int_{\R^n}|u|\,dx\leq 0
\qquad\text{in the sense of distributions in $(-\infty,T)$}.
$$
This proves that $u\equiv 0$. Notice also that a slight modifications of this argument, considering 
$\sqrt{(z^+)^2+\epsilon^2}-\epsilon$ instead of $\beta_\eps$, provides a comparison principle for the continuity
equation.

\begin{remark}[Special classes of vector fields] {\rm Many classes of vector fields are by now known to satisfy the
assumption of the previous theorem. Besides the case of Sobolev vector fields, treated in \cite{DiPerna-Lions89}, extensions
have been provided to $BV$ vector fields in \cite{Ambrosio03}. It is also interesting to note that well-posedness
results can be obtained for vector fields with a special structure, for instance $B(t,x)=(\bb(x),\nabla\bb(x)p)$ (with $\bb$
only Sobolev!) \cite{LeBris-Lions}. We refer to the Lecture Notes (in chronological order) \cite{cetraro}, \cite{bologna},
\cite{edinburgo} for much more on this topic.
}\end{remark}

So, we have seen that the key property is the renormalization property. So, the next question is: where this property
comes from, and how is it related to the regularity of $\bb_t$? To understand this, let us mollify not the vector field (as
we did for existence) but the weak solution $u_t$, and let us try to write down a PDE for $u^\eps_t$. We find
\begin{equation}\label{eq:5}
\frac{d}{dt}u^\eps_t+\bb_t\cdot\nabla u^\eps_t=-(\div\bb_tu_t)\ast\rho_\eps+\Commutator^\eps(\bb_t,u_t)
\end{equation}
where $\Commutator^\eps(\bb_t,u_t)$ is the so-called commutator between mollification and gradient along $\bb_t$
(notice that $\Commutator^\eps$ acts only on the times slices, i.e.\ $t$ by $t$):
$$
\Commutator^\eps(\cc,v):=\cc\cdot\nabla (v\ast\rho_\eps)-\bigl(\cc\cdot\nabla v\bigr)\ast\rho_\eps.
$$
Since $u^\eps_t$ are smooth, we can multiply both sides in \eqref{eq:5} by $\beta'(u_t^\eps)$ to get
$$
\frac{d}{dt}\beta(u^\eps_t)+\bb_t\cdot\nabla\beta(u^\eps_t)=-\beta'(u^\eps_t)(\div\bb_tu_t)\ast\rho_\eps+
\beta'(u^\eps_t)\Commutator^\eps(\bb_t,u_t).
$$
Hence, if we are able to prove that $\Commutator^\eps(\bb_t,u_t)\to 0$ strongly as $\eps\searrow 0$, we can take limits as $\eps\searrow 0$ in
the previous equation to obtain \eqref{eq:4}. 

In order to prove the strong convergence of the commutators, passing to the conservative form, we can also write
\begin{eqnarray*}
\tilde\Commutator^\eps(\cc,v)&:=&\div (\cc (v\ast\rho_\eps))-(\div (\cc v))\ast\rho_\eps\\&=&
\Commutator^\eps(\cc,v)+(v\ast\rho_\eps)\div\cc-(v\div\cc)\ast\rho_\eps,
\end{eqnarray*}
therefore the strong convergence as $\eps\searrow 0$ can be studied equivalently in the conservative
or non-conservative form. Now, if kernel $\rho$ is even, an integration by parts argument provides the representation
$$
\tilde\Commutator^\eps(\cc,v)(x)=\int_{\R^n}v(x-\eps y)\langle\frac{\cc(x+\eps y)-\cc(x)}{\eps},\nabla\rho(y)\rangle\,dy
+v\ast\rho_\eps(x)\div\cc(x).
$$
Using this representation it is not hard to show that the convergence is strong when
$\cc$ is Sobolev, because the strong convergence of difference quotients (here averaged w.r.t.\ to all directions $y$, with weight
given by $\nabla\rho$), characterizes Sobolev spaces. We know anyhow that the limit should be 0, but this can be seen
immediately because an integration by parts gives
$$
\lim_{\eps\searrow 0}\tilde\Commutator^\eps(\cc,v)(x)=v(x)\int_{\R^n}\langle\nabla\cc(x)y,\nabla\rho(y)\rangle\,dy+
v(x)\div\cc(x)=0.
$$ 
As this simple computation shows, cancellations can play an important role in the analysis of the asymptotic behaviour
of the commutators. This will be even more clear in the abstract setup, where an interpolation technique will lead to the
estimate of the commutators. This concludes the sketch of proof of Theorem~\ref{thm:maindp}. 

However, this approach is ``too Euclidean'' and not ``coordinate-free''. For instance, the theory can be immediately
adapted to compact Riemannian manifolds, working in local coordinates. But, if the manifold is unbounded, patching the
estimates in the different charts might be a problem. It is then natural to look for a more global and synthetic approach, suitable
for the extension also to nonsmooth and infinite-dimensional spaces. The key idea, already implemented in \cite{AF-09} (and then
in \cite{Trevisan-13}, for $BV$ vector fields), is
to use the heat semigroup as regularizing operator. However, while in \cite{AF-09,Trevisan-13} the computation were based on an explicit
knowledge of the transition probabilities (via Mehler's formula for the Ornstein-Uhlenbeck semigroup), in \cite{Ambrosio-Trevisan}
we obtain a new formula for the commutator, which provides strong convergence as soon as the semigroup has sufficiently
strong regularizing properties. 

We will show in the abstract part (see \eqref{eq:new_commu}) that, when the mollification is provided by the heat semigroup, 
in the divergence-free case the commutator can be expressed in integral form by
$$
\int_{\R^n}\Commutator^\eps(\cc,v)w\,dx=\frac 12
\int_0^\eps\biggl(\int_{\R^n} \nabla^{sym}\cc(\nabla P_{\eps-s}v,\nabla P_sw)\,dx\biggr)\,ds
$$
and we can use this new representation to estimate it.

\section{Abstract setup}

Our abstract setup in these notes will the standard one of the theory of Dirichlet forms: standard references are \cite{Bouleau-Hirsch91, Fukushima-Oshima-Takeda11}; we also point out the recent monograph \cite{BGL-13}, focused on $\Gamma$-calculus techniques and their applications to geometric and
functional inequalities. We have a Polish space $(X,\tau)$ and a Borel probability measure
$\mm$, with $\supp\mm=X$. Then we consider a $L^2$-lower semicontinuous quadratic form $\cE:L^2(X,\mm)\to [0,\infty]$, 
whose domain $D(\cE)$ is dense in $L^2(X,\mm)$. We still use the notation form $\cE(f,g)$ for the bilinear form associated to
$\cE$, defined in $D(\cE)\times D(\cE)$. In the sequel we denote
$$
\V:=D(\cE),\qquad \V_\infty:=D(\cE)\cap L^\infty(X,\mm).
$$
In addition we assume:

$\bullet$ Markovianity: $\cE(\eta\circ f)\leq\cE(f)$ for all $\eta:\R\to\R$ $1$-Lipschitz;

$\bullet$ Conservation of mass:  $\cE(1)=0$; 

$\bullet$ Strong locality: $\cE(f,g)=0$ if $fg=0$ $\mm$-a.e.\ in $X$;

$\bullet$ Existence of a Carr\'e du champ, namely a symmetric bilinear form $\Gamma:\V_\infty\times\V_\infty\to L^1(X,\mm)$
satisfying
$$
\cE(f,f\phi)-\frac 12\cE(f^2,\phi)=\int_X\Gamma(f)\phi\,d\mm.
$$
One can then canonically extend $\Gamma$ to a symmetric bilinear form on $\V\times\V$ such that
$\cE(f,g)=\int_X\Gamma(f,g)\,d\mm$ and prove the calculus rules
$$
\Gamma(\eta\circ f,g)=\eta'(f)\Gamma(f,g),\qquad
\Gamma(fg,h)=f\Gamma(g,h)+g\Gamma(f,h).
$$
Using the Carr\'e du champ we can also consider, for $p\in [2,\infty]$, the subspaces
\begin{equation}\label{def_Vp}
\V_p:=\left\{f\in \V:\ \Gamma(f)\in L^{p/2}(X,\mm)\right\}
\end{equation}
endowed with the norm $\|f\|_{\V_p}:=\|f\|_2+\|\Gamma(f)\|_{p/2}$.

\smallskip
Let us discuss now two important classes of examples.

\begin{example} [Infinitesimally Hilbertian metric measure spaces]\label{ex:infinitesimally-hilbertian}{\rm
Let $(X,d)$ be a complete and separable metric space and a Borel probability measure $\mm$ in $X$. For any
$f\in {\rm Lip}(X,d)$ we denote by $|\nabla f|$ the slope (also called local Lipschitz constant)
$$
|\nabla f|(x):=\limsup_{y\to x}\frac{|f(y)-f(x)|}{d(y,x)}.
$$
We call Cheeger energy the $L^2(X,\mm)$ relaxation of the previous energy, namely
$$
{\sf Ch}(f):=\inf\left\{\liminf_{n\to\infty}\int_X|\nabla f_n|^2\,d\mm:\ \text{$f_n\in {\rm Lip}(X)$, $f_n\to f$ in $L^2(X,\mm)$}\right\}.
$$
The functional ${\sf Ch}(f)$ can be represented by integration the square of a local object, called minimal relaxed slope and denoted
$|\nabla f|_w$, defined as follows: if $RS(f)$ denotes the class of functions $g\geq h$, with $h$ weak limit point in $L^2(X,\mm)$
of $|\nabla f_n|$ for some sequence $f_n\to f$ in $L^2(X,\mm)$, it is not hard to show (see \cite{AGS11a} for details)
that $RS(f)$ is a convex closed set in $L^2(X,\mm)$, not empty if and only if ${\sf Ch}(f)<\infty$. The minimal relaxed slope is precisely the element
with smallest norm in $RS(f)$, and 
$$
{\sf Ch}(f)=\int_X|\nabla f|_w^2\,d\mm.
$$

In general ${\sf Ch}$, though 2-homogeneous and convex, is not a quadratic form (for instance one can take
$X=\R^n$ and the distance induced by a non-Hilbertian norm, whose dual is non-Hilbertian as well). We call
a metric measure space $(X,d,\mm)$ \textit{infinitesimally Hilbertian} if ${\sf Ch}$ is a quadratic form.
The connection with the theory of Dirichlet forms has been provided in \cite{AGS11b}, where it is proved that,
for asymptotically Hilbertian metric measure spaces, ${\sf Ch}$ is a Markovian, strongly local and quasi regular
Dirichlet form, whose Carr\'e du champ is precisely $|\nabla f|_w^2$.
}\end{example}

\begin{example}[Log-concave measures]{\rm
Let $H$ be a separable Hilbert space and let $\gamma$ be a Borel probability measure on $H$. We say that 
$\gamma$ is log-concave if 
$$
\log \gamma\bigl((1-t)A+tB\bigr)\geq (1-t)\log\gamma(A)+t\log\gamma(B)
$$
for all pairs of open sets $A,\,B\subset H$. The class of log-concave measures has been widely studied, since the
work of Borell, see for instance \cite{Bogachev-10}, and it includes Gaussian measures $\mm$, as well their multiplicative perturbations of
the form $e^{-V}\mm$, with $V$ convex.

We consider the algebra $\Algebra$ of bounded, smooth cylindrical functions and the quadratic form
(where $|\nabla f|^2=\sum_i|\partial_i f|^2$ is independent of the chosen orthonormal basis)
$$
\cE_0(f):=
\begin{cases}
\int_H|\nabla f|^2\,d\gamma&\text{if $f\in\Algebra$;}
\\
+\infty & \text{if $f\in L^2(H,\gamma)\setminus\Algebra$.}
\end{cases}
$$
The largest $L^2(H,\gamma)$ envelope 
$$\cE(f):=\inf\left\{\liminf_{n\to\infty}\cE_0(f_n): \text{$f_n\to f$ in $L^2(H,\gamma)$}\right\}
$$
is easily seen to be a quadratic form and it is $L^2(H,\gamma)$ lower semicontinuous, by construction. It turns out that
the log-concavity assumption on $\gamma$ ensures that $\cE$ is an extension of $\cE_0$ (i.e.\ $\cE_0$ is \textit{closable})
and using this information it is also easy to prove existence of the Carr\'e du champ, which reduces to $|\nabla f|^2$
if $f\in\Algebra$.

It is interesting to remark that the traditional proofs of closability use integrations by parts formulae (see for instance \cite{Albeverio-Kusuoka})
along quasi-invariant directions, but the existence of these special directions (which span the Cameron-Martin space in
the Gaussian case) is still an open problem, see Section~4.3 of \cite{Bogachev-10}. Nevertheless, using optimal transportation tools, which yield
that $4\cE_0(\sqrt{\rho})$ is the descending slope of relative entropy functional $\rho\mapsto\int_X\rho\ln\rho\,d\gamma$ with respect to $\gamma$ (which,
by convexity of the relative entropy in the Wasserstein space $\Probabilities{H}$, is indeed lower semicontinuous) 
the closability of $\cE_0$ has been proved in \cite{AmbrosioSavareZambotti09} with no restriction on $\gamma$, besides log-concavity.  
}\end{example}

Continuing our presentation of the setup, we shall denote by $\Delta$ the infinitesimal generator of the semigroup $P_t$, with domain
$D(\Delta)$. The former
can be for instance defined by
$$
\int_X g\Delta f\,d\mm=-\cE(f,g)\qquad\forall g\in\V 
$$
and it satisfies the chain rule
$$
\Delta \eta\circ f=(\eta'\circ f)\Delta f+(\eta''\circ f)\Gamma(f)
$$
for all $\eta\in C^2(\R)$ with bounded first and second derivatives.
The semigroup can be defined by the property $P_tf\in H^1_{\rm loc}((0,\infty);L^2(X,\mm))$ with
$$
\frac{d}{dt}P_t f=\Delta P_t f\qquad\text{for a.e.\ $t>0$.}
$$
By the Markov property, it turns out that $P_t$ satisfies the maximum and minimum principle, $P_t1=1$ and $P_t$ is contractive in all $L^p\cap L^2(X,\mm)$ spaces, $1\leq p<\infty$, so we can canonically extend $P_t$ to a contractive semigroup
in all $L^p(X,\mm)$ spaces. Furthermore, the following regularizing effects play an important role:
$$
\cE(P_t f)\leq \inf_{v\in\V} \cE(v)+\frac {1}{2t}\|v-f\|_2^2\qquad\forall t>0,
$$
\begin{equation}\label{eq:regularization2}
\|\Delta P_t f\|_2^2\leq \inf_{v\in D(\Delta)}\|\Delta v\|_2^2 +\frac {1}{t^2}\|v-f\|_2^2\qquad\forall t>0.
\end{equation}
Finally, the theory of analytic semigroups provides also the estimate, for $p \in (1, \infty)$,
$$
\|\Delta P_t f\|_p\leq \frac{c^\Delta_p}{t}\|f\|_p\qquad\forall f\in L^p(X,\mm),\,\,t>0.
$$

Finally, we consider an algebra $\Algebra\subset\V_\infty$ of ``nice'' functions, dense in $\V$, stable under the joint composition
$\phi(f_1,\ldots,f_n)$ with $f_i\in\Algebra$ and $\phi\in C^1(\R^n)$ and Lipschitz (in particular $1\in\Algebra$). For instance
in infinite-dimensional Gaussian (or log-concave) spaces a typical choice of $\Algebra$ is provided by the smooth cylindrical functions,
while in infinitesimally Hilbertian metric measure spaces the natural choice of $\Algebra$ is provided by bounded Lipschitz functions. 

\section{Derivations and their regularity}

The concept of derivation in a measure space has been introduced and deeply studied in Chapter 7 of \cite{weaver1}, see also
\cite{weaver}. In the more recent papers \cite{Bate-12}, \cite{DiMarino-14}, \cite{gigli-WIP}, \cite{Schioppa-13}, derivations have 
been proved to be the natural tool for the development of a differential calculus in metric measure spaces, including the theory
of Sobolev spaces, differentiability of Lipschitz functions and second-order calculus.

\begin{definition}[Derivation] \label{def:derivation}We say that a linear functional $\bb:\Algebra\to L^0(X,\mm)$ is a derivation if
$\bb$ satisfies the Leibniz rule, namely
$$
\bb(fg)=f\bb(g)+g\bb(f)\qquad\forall f,\,g\in\Algebra.
$$
Given a derivation $\bb$, we write $|\bb|\in L^p(X,\mm)$ if there exists $g\in L^p(X,\mm)$ satisfying
$$
|\bb(f)|\leq g\sqrt{\Gamma(f)}\,\,\text{$\mm$-a.e.\ in $X$, for all $f\in\Algebra$.}
$$
We denote by $|\bb|$ the smallest function (in the $\mm$-a.e.\ sense) having this property.
\end{definition}

A typical example of derivation is provided by gradient vector fields, namely
$$
\bb_V(f)=\Gamma(V,f)\qquad f\in\Algebra
$$
with $V\in\V$. It is immediately seen that $|\bb_V|=\sqrt{\Gamma(V)}$, hence $|\bb_V|\in L^2(X,\mm)$. Since the space of
derivations has a natural structure of $L^\infty$ module, all finite combinations of the form
$$
\sum_{i=1}^N\chi_i\bb_{V_i}\qquad V_i\in\V,\,\,\chi_i\in L^\infty(X,\mm)
$$
are derivations.

\begin{definition} [Divergence] Let $\bb$ be a derivation with $|\bb|\in L^1(X,\mm)$. We say that $\div\bb$ belongs to
$L^p(X,\mm)$ if there exists $g\in L^p(X,\mm)$ satisfying
$$
\int_X \bb(f)\,d\mm=-\int_X gf\,d\mm\qquad\forall f\in\Algebra.
$$
By a density argument, it is easily seen that $g$, if exists, is unique.
The unique $g$ satisfying this property will be denoted $\div\bb$.
\end{definition} 

For gradient derivations we obviously have, by the very definition of $\Delta$,
$$
\div\bb_V=\Delta V.
$$

\begin{remark}[On the structure of the space of derivations]{\rm 
Using $|\bb|$ one can endow the vector space of derivations in $L^2$ with the structure of Hilbert module; using this structure, it
can be proved that the $L^\infty$ module generated by gradient derivations, namely
$$
\left\{\sum_{i=1}^N \chi_i\bb_{V_i}:\ \chi_i\in L^\infty(X,\mm),\,\,V_i\in\V\right\}
$$
is dense, see \cite{weaver} and \cite{gigli-WIP} for details.}
\end{remark}

As in the Euclidean theory, we need some regularity of the derivation in order to prove well-posedness
of the continuity equation. If we are in the category of Sobolev spaces, it was already understood in
\cite{Capuzzo} that only the symmetric part of the gradient plays a role. The next definition, which
is the natural extension of Bakry's weak definition of Hessian \cite{Bakry-94} from gradient vector fields to
general derivations, is consistent with the classical picture of the canonical Dirichlet form in a Riemannian manifold.

\begin{definition}[Derivations with deformation in $L^2$]
Let $\bb$ be a derivation in $L^2$, with $\div\bb\in\Lbm 2$, and assume that $\Algebra$ is dense in $\V_4$ (recall \eqref{def_Vp}). We write $D^{sym}\bb\in L^2(X,\mm)$ if
there exists $c \ge 0$ satisfying
\begin{equation}\label{eq:ineq_deformation}
\biggl|\int D^{sym}\bb(f,g) d\mm \biggr| \le c \|\sqrt{\Gamma(f)}\|_4 \|\sqrt{\Gamma(g)}\|_4  ,
\end{equation}
for all $f,\,g\in\V_4$ with $\Delta f,\,\Delta g\in L^4(X,\mm)$, where
	\begin{equation}
		\label{eq:distributional-deformation}
		 \int D^{sym}\bb(f,g) d\mm := -\frac{1}{2} \int \sqa{ \bb(f)\Delta g + \bb(g) \Delta f - (\div\bb) \Gamma(f,g) } d\mm.
	\end{equation}
	We let $\nor{D^{sym}\bb}_{2}$ be the smallest constant $c$ in \eqref{eq:ineq_deformation}.
\end{definition}

The density assumption on $\Algebra$ in the definition above allows for a unique continuous 
extension of $\bb$ to a derivation defined on $\V_4$, so that \eqref{eq:distributional-deformation} is well-defined.

\section{Eulerian side}

\subsection{Existence of solutions to the continuity equation}

In this section, we fix $T\in (0,\infty)$, a weakly$^*$-measurable family of derivations $\bb_t$, $t\in (0,T)$ (i.e.
$(t,x)\mapsto\bb_t(f)(x)$ is $\Leb{1}\times\BorelSets{\tau}$-measurable for all $f\in\Algebra$) and we consider
the continuity equation with a source term, namely 
\begin{equation}\label{eq:CEbis}
\frac{d}{dt}u_t+\div(\bb_tu_t)=c_t u_t.
\end{equation}
Under the minimal integrability assumptions 
$$
u\in L^\infty(L^p), \qquad |\bb|\in L^1(L^{p'}),\qquad c\in L^1(L^{p'})
$$
(always satisfied in the cases we shall consider), we can give
a meaning to \eqref{eq:CEbis} via its weak formulation, i.e.
$$
\frac{d}{dt}\int_X u_t f\,d\mm=\int_X (\bb_t(f)+c_tu_tf)\,d\mm\quad\text{in $\mathcal D'(0,T)$, for all $f\in\Algebra$.}
$$
Equivalently, we can use also test functions in the $t$ variable, writing
\begin{equation}\label{eq:6}
\int_0^T\int_X[-\psi'(t)u_tf+\psi(t)(\bb_t(f)+c_tu_tf)]\,d\mm dt=0\quad\text{for all $\psi\in C^1_c(0,T)$, $f\in\Algebra$.}
\end{equation}
Notice that, in these weak formulations, $u_t$ is only uniquely determined up to a $\Leb{1}$-negligible set of times. However, using
the fact that $t\mapsto\int u_tf\,d\mm$ has a continuous representative for all $f\in\Algebra$, it is not hard to prove
the existence of a unique representative $t\mapsto u_t$ which is continuous in the duality with $\Algebra$, and we shall always
tacitly refer to this distinguished representative.

We can encode an initial condition $\bar u$ by requiring, for instance, that $\int_X u_t f\,d\mm\to\int_X \bar u f\,d\mm$
as $t\downarrow 0$ for all $f\in\Algebra$. For technical reasons, we will also use a variant of \eqref{eq:6}, writing
\begin{equation*}
\int_0^T\int_X[-\psi'(t)u_tf+\psi(t)(\bb_t(f)+c_tu_tf)]\,d\mm dt=\psi(0)\int_Xf\bar u\,d\mm
\quad\text{$\forall\psi\in C^1([0,T])$, $\psi(T)=0$}
\end{equation*}
for all $f\in\Algebra$.

Notice that the condition $1\in\Algebra$ ensures that any weak solution $u$ satisfies, when $c_t=0$, the mass conservation property
$$
\int_X u_t\,d\mm=\int_X\bar u\,d\mm\qquad\forall t\in [0,T].
$$

\begin{theorem} \label{thm:existence} Let $\bar u\in L^r(X,\mm)$ for some $r\in [2,\infty]$. Assume that $|\bb|\in L^1(L^2)$,
$\div\bb\in L^1(L^2)$ and that its negative part $(\div\bb)^-$ belongs to $L^1(L^\infty)$. Then there
exists a solution $u\in L^\infty(L^r)$ to \eqref{eq:CEbis} with $c_t=0$. In addition, we can build $u$ in such a way that
$u\geq 0$ if $\bar u\geq 0$.
\end{theorem}

The proof of the theorem can be obtained by the classical method of vanishing viscosity in three steps: first we obtain a solution to a regularized
equation, then we prove a priori estimates in $L^r$, eventually we take the limit as the viscosity parameter tends to 0. We sketch the
proof under the simplifying additional assumption that $|\bb|\in L^\infty(L^2)$ and $(\div\bb)^-\in L^\infty(L^2)$.

\noindent
{\bf Step 1.} (Existence for the regularized equation) For $\sigma\in (0,1)$, let us prove existence of a solution $u^\sigma$ to
\begin{equation}\label{eq:viscous}
\frac{d}{dt}u_t+\div(\bb_t u_t)=\sigma u_t,\qquad u_0=\bar u.
\end{equation}
We use the well-known J.L.Lions' lemma, a generalization of Lax-Milgram theorem, whose statement is recalled below (see for instance 
\cite[Thm.~III.2.1, Corollary III.2.3]{Showalter97} for the proof).

\begin{lemma}[Lions]\label{lem:lions}
Let $V$, $H$ be respectively a normed and a Hilbert space, with $V$ continuously embedded in $H$, with $\|v\|_H\leq\|v\|_V$
for all $v\in V$, and let
$B:V\times H\to\R$ be bilinear, with $B(v,\cdot)$ continuous for all $v\in V$. If $B$ is coercive, namely there exists $c>0$ satisfying
$B(v,v)\geq c\|v\|^2_V$ for all $v\in V$, then for all $\ell\in V'$ there exists $h\in H$ such that $B(\cdot,h)=\ell$ and
\begin{equation*}
\nor{h}_H\leq\frac{\nor{\ell}_{V'}}{c}.
\end{equation*}
\end{lemma}

We apply the lemma with $H=L^2((0,T);\V)=L^2(\V)$, with $V\subset H$ equal to the vector space generated by the functions $\psi f$, with
$\psi\in C^1_c([0,T))$ and $f\in\Algebra$, endowed with the norm
\begin{equation*}
	 \nor{\vphi}_V^2 = \|\varphi\|_{L^2(\V)}^2 + \|\varphi_0\|^2_2,
\end{equation*}
with the bilinear form
$$
B(\vphi,h) = \int_0^T\int_X \sqa{-\partial_t \vphi +\lambda \vphi -d\vphi(\bb)} h + \sigma \Gamma(\vphi,h) d\mm dt
$$
and eventually with $\ell(\varphi) = \int \vphi_0\bar u d\mm$. It is immediate to check, by the definition of the norm in
$V$, that $\ell\in V'$ and that $\|\ell\|_{V'}\leq\|\bar u\|_2$. Also the continuity of $B(\vphi,\cdot)$ is immediate to check, using the
growth assumption $|\bb|\in L^\infty(L^2)$.
Let us prove that $B$ is coercive for $\lambda$ sufficiently large. Indeed,
\begin{equation*}
	\begin{split}
		\int_0^T\int_X \sqa{-\partial_t\varphi+\lambda\varphi -\bb(\varphi)} \varphi\, d\mm dt
		& = \|\varphi_0\|_2^2+\lambda \nor{\varphi}^2_{L^2(L^2)} -\frac{1}{2}\int \bb(\varphi^2)\,d\mm dt\\
		& \ge \|\varphi_0\|_2^2+\lambda\nor{\varphi}^2_{L^2(L^2)} -\frac{1}{2}\int\varphi^2\dbneg d\mm dt \\
		& \ge \|\varphi_0\|_2^2+(\lambda-\frac 12 \|\dbneg\|_{\infty})\nor{\varphi}_{L^2(L^2)}^2\\
		& \geq\sigma (\|\varphi_0\|_2^2+\nor{\varphi}_{L^2(L^2)}^2)
	\end{split}
	\end{equation*}
	if $\lambda\geq\lambda_\sigma=\frac 12 \|\dbneg\|_{\infty}+\sigma$. Taking also into account the term $\sigma\int\Gamma(\varphi,\varphi)$
	appearing in $B(\vphi,\vphi)$ we obtain $B(\varphi,\varphi)\geq\sigma\|\vphi\|_V^2$. Hence,
By applying Lemma~\ref{lem:lions} with $\lambda=\lambda_\sigma$ we obtain a weak solution $h$ to
$$
\frac{d}{dt}h_t+\div(\bb_t h_t)+\lambda_\sigma h_t=\sigma h_t,\qquad h_0=\bar u
$$
with $\|h\|_H\leq \|\bar u\|_2/\sigma$. Hence, setting $u=e^{\lambda_\sigma t}h$, we obtain a solution to \eqref{eq:viscous}
satisfying
\begin{equation}\label{eq:first_a_priori}
\|e^{-\lambda_\sigma t} u\|_{L^2(\V)} \leq \frac{\|\bar u\|_2}{\sigma}\quad\text{with}\quad\lambda_\sigma=\frac 12\|(\div\bb_t)^-\|_\infty+\sigma.
\end{equation}

\noindent
{\bf Step 2.} (A priori estimates) The approximate solutions $u^\sigma$ of step 1 satisfy the apriori estimate
$$
\sup_{t\in (0,T)}\|(u_t^\sigma)^\pm\|_r\leq \|\bar u^\pm\|_r\exp\bigl((1-\frac 1r)\|(\div\bb_t)^-\|_{L^1(L^\infty)}\bigr).
$$
This can be achieved by a further approximation, i.e.\ by approximating $u^\sigma_t$ by $P_su^\sigma_t$. The strategy
is to differentiate w.r.t.\ $t$ the entropy $\int \beta(P_su^\sigma_t)\,d\mm$, letting first $\beta(z)\to |z^\pm|^r$ and then
$s\searrow 0$. Notice that these approximation provide also a kind of renormalization property; however,
the $L^2(\V)$ regularity of $u^\sigma$ makes the proof much easier and allows to dispense with
the regularizing assumption on the semigroup needed in the next section, where the solutions
have no regularity but satisfy only a uniform $L^r$ bound.

\noindent
{\bf Step 3.} (Limit as $\sigma\searrow 0$) In order to pass to the limit as $\sigma\searrow 0$ in \eqref{eq:viscous}
we use the weak formulation
$$
\frac{d}{dt}\int_X fu^\sigma_t\,d\mm=\int_X \bb_t(f)u^\sigma_t\,d\mm+\sigma\int_X\Gamma(u^\sigma_t,f)\,d\mm
\qquad \quad f\in\Algebra
$$
and the a-priori bound provided by \eqref{eq:first_a_priori},
which ensures that $\sigma u^\sigma$ is bounded in $L^2(\V)$. Since $\sigma u^\sigma$ strongly converge to 0 in $L^2(L^2)$,
we obtain that $\sigma u^\sigma$ weakly converge to $0$ in $L^2(\V)$. Therefore 
$$
\lim_{\sigma\searrow 0}\sigma\int_0^T\int_X\Gamma(u^\sigma_t,f)\,d\mm dt=0
$$
and this proves that any weak$^*$ limit as $\sigma\searrow 0$ of $u^\sigma$ in $L^\infty(L^r)$ is a solution to \eqref{eq:viscous}.
In an analogous way one can see that the also the initial condition $\bar u$, the sign condition (if any)
and the uniform $L^r$ bounds are retained in the limit.

\subsection{Uniqueness of solutions to the continuity equation}

For the results of this section we need a mild enforcement of the regularizing properties of the semigroup,
given in the next definition.

\begin{definition}[$L^p$-$\Gamma$ inequality]
Let $p\in [1,\infty]$. We say that the $L^p$-$\Gamma$ inequality holds if there exists $c_p>0$ satisfying
$$
\Big\|\sqGq{ P_t f} \Big\|_p \le \frac{c_p}{\sqrt{t}} \nor{f}_p, \quad\text{for every $f \in L^2\cap L^p(X,\mm)$,  $t \in (0,1)$.}
$$
\end{definition}

Although the $L^p$-$\Gamma$ inequality is expressed for $t \in (0,1)$, from its validity and $L^p$ contractivity of $P$, we easily deduce that
\begin{equation*} 
\nor{ \sqGq{ P_t f} }_p \le \frac{c_p }{\sqrt{t}\land 1}  \nor{f}_p,\quad \text{for every $f \in L^2\cap L^p(X,\mm)$, $t \in (0,\infty)$.}
\end{equation*}
Notice also that, thanks to \eqref{eq:regularization2}, the $L^2$-$\Gamma$ inequality always holds, with $c_2=1/\sqrt{2}$. 
By semilinear Marcinkiewicz interpolation, 
we obtain that if the $L^p$-$\Gamma$ inequality holds then, for every $q$ between $2$ and $p$, the $L^q$-$\Gamma$ inequality holds as well.

\begin{theorem} \label{thm:well_PDE} Assume that $\Algebra$ is dense in $\V_4$ and that $L^4$-$\Gamma$ inequalities hold. If
$$
|D^{sym}\bb_t|\in L^1(L^2),\qquad
|\bb|\in L^1(L^2),\qquad
|\div\bb_t|\in L^1(L^\infty)
$$
then  the continuity equation \eqref{eq:CE} with initial datum $\bar u\in L^2(X,\mm)$ is well posed in $L^\infty(L^4)$.
\end{theorem}

We already know, from the previous section with $r=4$, about existence.
As in the Euclidean case, in order to prove uniqueness we need to prove that 
weak solutions to \eqref{eq:CE} are renormalized, i.e.
\begin{equation}\label{eq:5ter}
\frac{d}{dt}\beta(u_t)+\bb_t\cdot\nabla \beta(u_t)=-u_t\beta'(u_t)\div\bb_t\qquad\forall\beta\in C^1(\R)\cap {\rm Lip}(\R).
\end{equation}
Still following the Euclidean case, we prove this property by regularization, i.e.\ writing $u^\alpha_t=P_\alpha u_t$ and looking at
the PDE satisfied by $t\mapsto u^\alpha_t$, namely
\begin{equation}\label{eq:5bis}
\frac{d}{dt}u^\alpha_t+\bb_t\cdot\nabla u^\alpha_t=-u^\alpha_t\div\bb_t+\Commutator^\alpha(\bb_t,u_t)
\end{equation}
where $\Commutator^\alpha(\cc,v)$ is defined by 
$$
\Commutator^\alpha(\cc,v):=\div \bigl((P_\alpha v)\cc\bigr)-P_\alpha^*(\div(v\cc)\bigr).
$$
We used the notation $P_\alpha^*$ because, strictly speaking, $\div (v\cc)$ is not a function, so the action of the
semigroup should be understood in the dual sense (this will also be clear from the proof of the commutator estimate).
If we are able to show that
\begin{equation*}
\lim_{\alpha\searrow 0}\int_0^T\|\Commutator^\alpha(u_t,\bb_t)\|_1\,dt=0
\end{equation*}
then we can pass to the limit as $\alpha\searrow 0$ in 
$$
\frac{d}{dt}\beta(u^\alpha_t)+\bb_t\cdot\nabla\beta(u^\alpha_t)=-u^\alpha_t\beta'(u^\alpha_t)\div\bb_t
+\beta'(u^\alpha_t)\Commutator^\alpha(\bb_t,u_t),
$$
derived from \eqref{eq:5bis} thanks to the higher regularity of $u^\alpha_t$, to obtain \eqref{eq:5ter}.

The proof of \eqref{eq:5bis} would be trivial if we knew that $u_t\in\Algebra$; indeed, for $v\in\Algebra$ one has
$$
\Commutator^\alpha(\cc,v)=(P_\alpha v)\div\cc+\cc(P_\alpha v)-P_\alpha(v\div\cc)-P_\alpha (\cc(v))
$$
and it is easily seen that this expression converges to $0$ in $L^1(X,\mm)$ as $\alpha\searrow 0$.
In general, since we only know that $u_t\in L^2(X,\mm)$, \eqref{eq:5bis}, is a consequence of a density argument and of the following commutator estimate.

\begin{theorem} [Commutator estimate]\label{thm:commutator} Assume that $\Algebra$ is dense in $\V_4$ and that $L^4$-$\Gamma$ inequalities hold. 
Then, there exists a constant $c$ satisfying
\begin{equation}\label{eq:stima_commutatori}
\|\Commutator^\alpha(\cc,v)\|_{4/3}\leq c\|v\|_4\bigl[\|D^{sym}\cc\|_2+\|\div\cc\|_2\bigr] \quad \forall \alpha \in (0,1),
\end{equation}
for all derivations $\cc$ with $|\cc|\in L^2(X,\mm)$ and $|D^{sym}\cc|\in L^2(X,\mm)$ and all $v\in L^4(X,\mm)$.
\end{theorem}

The proof of the commutator estimate relies on Bakry's interpolation, very much as in the derivation of gradient
contractivity $\Gamma(P_tf)\leq e^{-2Kt}P_t\Gamma(f)$ from Bochner's inequality
$$
\Delta\frac 12\Gamma(f)-\Gamma(f,\Delta f)\geq K\Gamma(f).
$$
We provide the proof in the simplified case when $\div\cc=0$ and $v\in\Algebra$, assuming for simplicity that
$\Algebra$ is invariant under the action of the semigroup, as it happens in many cases of interest. 
We refer to \cite{Ambrosio-Trevisan} for the complete proof, where the assumption that
$\Algebra$ is dense in $\V_4$ plays a role (as in the Euclidean case, when $\div\cc\in L^\infty(X,\mm)$
lower order terms appear in the computation below). Let us write
\begin{eqnarray*}
\Commutator^\alpha(\cc,v)&=&\int_0^\alpha\frac{d}{ds}P_{\alpha-s}^*\bigl(\div(P_s v\cc)\bigr)\,ds\\&=&
\int_0^\alpha -\Delta^* P_{\alpha-s}\bigl(\div(P_s v\cc)\bigr)+P_{\alpha-s}^*\bigl(\div (\Delta P_s v\cc)\bigr)\, ds.
\end{eqnarray*}
Now, instead of expanding the divergence terms, we take $w\in\Algebra$ and use the dual form of the
commutator to get
\begin{eqnarray*}
\langle \Commutator^\alpha(\cc,v),w\rangle&=&
\int_0^\alpha -\div(P_s v\cc)\Delta P_{\alpha-s}w+\div (\Delta P_s v\cc)P_{\alpha-s}w\, ds\\&=&
-\int_0^\alpha \cc(P_s v)\Delta P_{\alpha-s}w+\cc(P_{\alpha-s}w)P_s\Delta v\, ds.
\end{eqnarray*}
Taking \eqref{eq:distributional-deformation} into account, we get
\begin{equation}\label{eq:new_commu}
\langle \Commutator^\alpha(\cc,v),w\rangle=\frac 12\int_0^\alpha D^{sym}\cc(P_s v,P_{\alpha-s}w)\,ds
\end{equation}
and therefore the $L^4$ $\Gamma$-inequality combined with the assumption $\|D^{sym}\cc\|\in L^2(X,\mm)$ give
the estimate 
\begin{eqnarray*}
|\langle \Commutator^\alpha(\cc,v),w\rangle|&\leq&\frac 12\|D^{sym}\cc\|_2\int_0^\alpha\|\sqrt{\Gamma(P_sv)}\|_4
\|\sqrt{\Gamma(P_{\alpha -s}w)}\|_4\,ds\\&\leq&
c\|D^{sym}\cc\|_2\|v\|_4\|w\|_4\int_0^\alpha \frac{1}{\sqrt{s(\alpha-s)}}\,ds.
\end{eqnarray*}
Since $\Algebra$ is dense in $L^4(X,\mm)$, by duality this proves \eqref{eq:stima_commutatori}.

\section{Lagrangian side}

This section of the notes is devoted to the Lagrangian side of the theory. First we see how the definitions of ODE and
regular Lagrangian flow can be immediately adapted to the abstract setting. Then, in the subsequent two subsections,
we illustrate the two basic mechanisms which allow to move from the Lagrangian to the Eulerian side, and conversely.
Eventually these results will be applied to show existence and uniqueness of regular Lagrangian flows. 

\subsection{ODE's associated to derivations and Regular Lagrangian Flows}\label{sec:ODE-RLF}

In the sequel it will be useful to fix a countable subset $\Algebra^*\subset C(X)\cap\{f\in\Algebra:\ \Gamma(f)\leq 1\}$,
assuming that:
\begin{itemize}
\item[(a)] $\Q\Algebra^*$ is dense in $\V$ and a vector space over $\Q$;
\item[(b)] The possibly infinite distance
(following \cite{Biroli-Mosco95}, see also \cite{Sturm95,Stollmann10})
$$
d_{\Algebra^*} (x,y):=\sup_{f\in\Algebra^*}|f(x)-f(y)|
$$
satisfies $\lim\limits_{m,n\to\infty} d_{\Algebra^*}(x_n,x_m)\to 0$ implies $x_n\to x$ in $(X,\tau)$ for some $x\in X$.
\end{itemize}
We adopt the convention that, whenever $f\in\Algebra^*$, $f$ stands for the unique continuous
representative. 

For instance, in the metric measure setup one can choose as generators of $\Algebra^*$ the distance functions
from a countable set of points; it is easily seen in this case that $d_{\Algebra^*}$ coincides with $d$. Notice also
that, since we are assuming already that $(X,\tau)$ is complete, the condition in (b) can be stated
by saying that Cauchy sequences w.r.t.\ $d_{\Algebra^*}$ are Cauchy w.r.t.\ $\tau$, and a necessary
condition for this to hold is that $\Algebra^*$ separates points. 

\begin{definition}[ODE induced by a family of derivations]\label{def:ourODE}
Let $\eeta\in\Probabilities{C([0,T];(X,\tau)}$ and let $(\bb_t)_{t\in (0,T)}$ be a Borel family of derivations. We say that
$\eeta$ is concentrated on solutions to the ODE $\dot\eta=\bb_t(\eta)$ if 
 $$
 \text{$f\circ\eta\in W^{1,1}(0,T)$ and $\frac {d}{dt} (f\circ\eta)(t)=\bb_t(f)(\eta(t))$, 
for a.e.\ $t\in (0,T)$,}
$$
for $\eeta$-a.e.\ $\eta \in C([0,T];(X,\tau_{\Algebra^*}))$, for all $f\in\Algebra$.
\end{definition}

Notice that, for $f\in\Algebra^*$, $f\circ\eta$ is continuous for $\eeta$-a.e.\ $\eta$. Hence the $W^{1,1}(0,T)$
improves to $AC([0,T])$ regularity when $f\in\Algebra^*$. it is easily seen that the dual norm $|\bb|_*$, namely the smallest 
function $g_*$ satisfying
$$
|\bb(f)|\leq g_*\quad\text{$\mm$-a.e.\ in $X$ for all $f\in\Algebra^*$}
$$
is smaller than $|\bb|$ whenever $|\bb|\in L^2(X,\mm)$. Under suitable regularity assumptions on the metric measure structure (see
\cite{Ambrosio-Trevisan}) it can be proved that it coincides with $|\bb|$.

Notice also that the property of being concentrated on solutions to the ODE implicitly depends on the choice of Borel representatives of
the maps $f$ and $(t,x)\mapsto\bb_t(f)(x)$, $f\in\Algebra$. As such, it should be handled with care. However, as we have seen
in the Euclidean case, Fubini's theorem ensures that the for the class of regular flows of Definition~\ref{def:dregflow} 
this sensitivity to the choice of Borel representatives disappears.

\begin{definition} [Regular flows]\label{def:dregflow} We say that $\XX: [0,T]\times X\to X$ is a regular flow (relative
to $\bb$) if the following two properties hold: 
\begin{itemize}
\item[(i)] $\XX(0,x)=x$ and $\XX(\cdot,x)\in C([0,T];(X,\tau))$ for all $x\in X$;
\item[(ii)] for all $f\in\Algebra$, $f(\XX(\cdot,x))\in W^{1,1}(0,T)$ and $\frac{d}{dt} f(\XX(t,x))=\bb_t(f)(\XX(t,x))$ 
for a.e.\ $t\in (0,T)$, for $\mm$-a.e.\ $x\in X$;
\item[(iii)] there exists a constant $C=C(\XX)$ satisfying $\XX(t,\cdot)_\#\mm\leq C\mm$ for all $t\in [0,T]$.
\end{itemize}
\end{definition}

The following simple lemma is the first basic transfer mechanism from the Lagrangian to the Eulerian side. The converse
mechanism, the so-called superposition principle, will be introduced in the next section.
The lemma shows that time marginals of measures $\eeta$ concentrated on solutions
to the ODE $\dot\eta=\bb_t(\eta)$ provide weakly continuous solutions to the continuity equation.

\begin{lemma}\label{eq:fromODEtoPDE}
Let $\eeta\in\Probabilities{C([0,T];(X,\tau))}$ be concentrated on solutions $\eta$ to the ODE $\dot\eta=\bb_t(\eta)$, where 
$|\bb|\in L^1(L^p)$ for some $p\in [1,\infty]$ and $\mu_t:=(e_t)_\#\eeta\in\Probabilities{X}$ are representable as 
$u_t\mm$ with $u \in L^\infty(L^{p'})$. Then, the following two properties hold:
\begin{itemize}
\item[(a)] the family $(u_t)_{t\in (0,T)}$ is a weakly continuous, in the duality with $\Algebra$, solution to the continuity equation;
\item[(b)] $\eeta$ is concentrated on $AC([0,T]; (X, d_{\Algebra^*}) )$, with
\begin{equation} \label{eq:inequality-or-equality}\abs{\dot \eta} (t) = \abs{\bb_t}_*(\eta(t))\quad \text{ for a.e.\ $t \in (0,T)$, for $\eeta$-a.e.\ $\eta$.} \end{equation}
\end{itemize}
\end{lemma}
\begin{proof} We integrate w.r.t.\ $\eeta$ the weak formulation
$$
\int_0^t -\psi'(t)f\circ\eta(t) \,dt=\int_0^T\psi(t)\bb_t(f)(\eta(t))\, dt
$$
with $f\in\Algebra$, $\psi\in C^1_c(0,T)$, to recover the weak formulation of the continuity equation for $(u_t)$.

Given $f \in \Algebra^*$, for $\eeta$-a.e.\ $\eta$, the map $t \mapsto f \circ \eta(t)$ is absolutely continuous, with
\[ f\circ\eta(t) - f\circ\eta(s)  = \int_s^t  \bb_r(f) (\eta(r)) dr, \quad\text{for all $s,\, t \in [0,T]$.}\]
By integration w.r.t.\ $\eeta$ we obtain that $t\mapsto\int_X fu_t\,d\mm$ is absolutely continuous in $[0,T]$ for all
$f\in\Algebra^*$, and a density argument gives the weak continuity in the duality with $\Algebra$.
In particular one has $\bb_t(f)(\eta(t))=(f\circ\eta)'(t)$ a.e.\ in $(0,T)$, for $\eeta$-a.e.\ $\eta$.

By Fubini's theorem and the fact that the marginals of $\eeta$
are absolutely continuous w.r.t.\ $\mm$ we obtain that, for $\eeta$-a.e.\ $\eta$, one has
\[ \sup_{f \in \Algebra^*} \abs{(f \circ \eta)'(t)} =  \sup_{f \in \Algebra^*} \abs{\bb_t(f)(\eta(t))} = \abs{\bb_t}_*(\eta(t)), \quad \text{for a.e.\ $t\in (0,T)$,}\]
and therefore
\[ d_{\Algebra^*} (\eta(t), \eta(s) ) = \sup_{f \in \Algebra^*}\bigl|(f\circ\eta)(t) - (f\circ\eta)(s)\bigr|
\le \int_s^t \abs{\bb_t}_*(\eta(r))dr, \quad \text{for all $s,\,t \in [0,T]$,} \]
proving that $\eta \in AC([0,T]; (X, d_{\Algebra^*}))$, with $\abs{\dot \eta}(t) \le \abs{\bb_t}_*(\eta(t))$, for a.e.\ $t\in (0,T)$. 
The converse inequality follows from the fact that every $f\in\Algebra^*$ is $1$-Lipschitz with respect to $d_{\Algebra^*}$, 
thus for $\eeta$-a.e.\ $\eta$ one has
\[ \abs{\bb_t}_*(\eta(t)) =  \sup_{f \in \Algebra^*} \abs{\bb_t(f)(\eta(t))} =
\sup_{f \in \Algebra^*} \abs{(f \circ \eta)'(t)} \le \abs{\dot \eta}(t), \quad \text{for a.e.\ $t\in (0,T)$.}\]
\end{proof}

\subsection{The superposition principles}

The superposition principle is one of the basic principles that allow to move from an Eulerian perspective
to a Lagrangian one. The principle goes back to L.C.Young and it has been implemented in various forms,
for instance in the context of the theory of currents \cite{Smirnov}. Stochastic counterparts of this principle are of
course the representation of families of random variables via stochastic processes, as for instance in the theory of Dirichlet
forms and Markov processes. For solutions to the continuity equation in $\R^n$ (somehow
a special class of $(n+1)$-dimensional space-time currents) the superposition principle can be stated in a particularly convenient form (where
also singular measures are allowed) see \cite[Thm.~8.2.1]{Ambrosio-Gigli-Savare05}
for the case when $|\bb|\in L^r(\nu_t dt)<\infty$ for some $r>1$, and \cite[Thm.~12]{bologna} for the case $r=1$.

\begin{theorem}[Superposition principle in $\R^n$]\label{thm:superpoRn}
Let $\bb:(0,T)\times\R^n\to\R^n$ be Borel and let
$\mu_t\in\Probabilities{\R^n}$, $t\in [0,T]$, be a weakly continuous, in the duality with 
$C_c(\R^n)$, weak solution to the continuity equation
$$
\frac{d}{dt}\mu_t+\div (\bb_t\mu_t)=0,
$$
with $\int_0^T\int_{\R^n}|\bb_t|\,d\mu_tdt<\infty$.
Then, there exists a Borel probability measure
$\llambda$ in $C([0,T];\R^n)$ satisfying $(e_t)_\#\llambda=\nu_t$ for all $t\in [0,T]$, concentrated on $\gamma\in AC([0,T];\R^n)$
which are solutions to the ODE $\dot\gamma=\bb_t(\gamma)$ a.e.\ in $(0,T)$. 
\end{theorem}

The strategy of proof is to mollify with a Gaussian kernel, getting $\mu^\eps_t=\rho^\eps_t\Leb{n}$ satisfying
$\rho^\eps_t>0$ and
$$
\frac{d}{dt}\mu^\eps_t+\div (\bb^\eps_t\mu^\eps_t)=0
$$
with $\bb^\eps_t:=(\bb_t\mu_t)\ast\rho_\eps/\mu_t\ast\rho_\eps$. Even though the Lipschitz bounds on $\rho^\eps_t$
are only local, it is possible to show that there exists a smooth flow $X^\eps$ associated to $\bb^\eps_t$ with
$X^\eps(t,\cdot)_\#\mu^\eps_0=\mu^\eps_t$. Denoting by $\llambda^\eps$ the probability measure image
of $\mu^\eps_0$ under the map $x\mapsto\XX^\eps(\cdot,x)$, the most delicate part of the proof is the limit as $\eps\to 0$:
since no continuity assumption is made on $\bb$, it is not obvious that weak limit points $\llambda$ are concentrated
on solutions to the ODE.

We endow $\R^\infty=\R^{\N}$ with the product topology and
we shall denote by $\pi^n:=(p_1,\ldots,p_n):\R^\infty\to\R^n$ 
the canonical projections from $\R^\infty$ to $\R^n$. On the space $\R^\infty$ 
we consider the complete and separable distance
$$
d_\infty(x,y):=\sum_{n=1}^\infty 2^{-n}\min \cur{1,|p_n(x)-p_n(y)|}.
$$
Accordingly, we consider the space $C([0,T];\R^\infty)$ endowed with the distance
$$
\delta(\eta,\tilde \eta):=\sum_{n=1}^\infty 2^{-n}\max_{t \in [0,T]}\min\cur{1,|p_n(\eta(t))-p_n(\tilde \eta(t))|},
$$
which makes
$C([0,T];\R^\infty)$ complete and separable as well.
We shall also consider the subspace $AC_w([0,T];\R^\infty)$ of $C([0,T];\R^\infty)$ consisting of all $\eta$ such that
$p_i\circ\eta\in AC([0,T])$ for all $i\geq 1$. Notice that for this class of curves the derivative $\eta'\in\R^\infty$ can still be
defined a.e.\ in $(0,T)$, arguing componentwise (we use the notation $AC_w$ to avoid the confusion with the space
of absolutely continuous maps from $[0,T]$ to $(\R^\infty,d_\infty)$).

We call regular cylindrical function any $f:\R^\infty\to\R$ representable in the form 
$$
f(x)=\psi(\pi_n(x))=\psi\bigl(p_1(x),\ldots,p_n(x)\bigr)\qquad x\in\R^\infty,
$$
with $\psi:\R^n\to\R$ bounded and continuously differentiable, with bounded derivative. Given $f$ regular cylindrical as above,
we define $\nabla f:\R^\infty\to c_0$ (where $c_0$ is the space of sequences $(x_n)$ null for $n$ large enough)
by
\begin{equation*}
\nabla f(x):=\bigl(\frac{\partial\psi}{\partial z_1}(\pi_n(x)),\ldots,\frac{\partial\psi}{\partial z_n}(\pi_n(x)),0,0,\ldots).
\end{equation*}

We fix a 
Borel vector field $\cc:(0,T)\times\R^\infty\to\R^\infty$ and a weakly continuous (in duality with regular cylindrical functions)
family of Borel probability measures $\{\nu_t\}_{t\in [0,T]}$ in $\R^\infty$ satisfying 
\begin{equation}\label{eq:conti1}
\int_0^T\int_{\R^\infty} |p_i(\cc_t)| \,d\nu_tdt <\infty,\qquad\forall i\geq 1
\end{equation}
and, in the sense of distributions,
\begin{equation}\label{eq:conti2}
\frac{d}{dt}\int_{\R^\infty} f\,d\nu_t=\int_{\R^\infty}\bra{\cc_t,\nabla f} \,d\nu_t\qquad\text{in $[0,T]$, for all $f$ regular cylindrical.}
\end{equation}

By a canonical cylindrical approximation combined with a tightness argument, Theorem~\ref{thm:superpoRn}
can be extended to $\R^\infty$ as follows.

\begin{theorem}[Superposition principle in $\R^\infty$]\label{thm:superpoRinfty}
Under assumptions \eqref{eq:conti1} and \eqref{eq:conti2}, there exists a Borel probability measure
$\llambda$ in $C([0,T];\R^\infty)$ satisfying $(e_t)_\#\llambda=\nu_t$ for all $t\in [0,T]$, concentrated on $\gamma\in AC_w([0,T];\R^\infty)$
which are solutions
to the ODE $\dot\gamma=\cc_t(\gamma)$ a.e.\ in $(0,T)$. 
\end{theorem}

Even though the distance $d_{\Algebra^*}$ may be equal to $\infty$, we call the next result ``superposition in metric measure spaces'', because in most cases $\Algebra^*$ consists precisely, as we already said, of distance functions from a countable 
dense set  (see also the recent papers \cite{Bate-12} and \cite{Schioppa-13} for related results on the existence of suitable measures in the space of curves, and derivations) and $d_{\Algebra^*}$ coincides with the original distance.

\begin{theorem}[Superposition principle in metric measure spaces]\label{thm:superpo}
Assume conditions (a), (b) above on $\Algebra^*$. 
Let $\bb=(\bb_t)_{t\in (0,T)}$ be a Borel family of derivations and 
let $\mu_t=u_t\mm\in\Probabilities{X}$, $0\leq t\leq T$, be a weakly continuous, in the duality
with $\Algebra$, solution to the continuity equation
\begin{equation*}
	\partial_t\mu_t + \div\bra{\bb_t\mu_t} =0
	\end{equation*} 
with
\begin{equation*}
u\in L^\infty_t(L^p_x),\qquad
\int_0^T\int|\bb_t|^rd\mu_tdt<\infty,\qquad \frac{1}{r}+\frac{1}{p}\leq 1/2.
\end{equation*}
Then there exists $\eeta\in\Probabilities{C([0,T];(X,\tau))}$ satisfying:
\begin{itemize}
\item[(a)] $\eeta$ is concentrated on solutions $\eta$ to the ODE $\dot\eta=\bb_t(\eta)$, 
according to Definition~\ref{def:ourODE};
\item[(b)] $\mu_t=(e_t)_\#\eeta$ for any $t\in [0,T]$.
\end{itemize}
\end{theorem}
\begin{proof} We enumerate by $f_i$, $i\geq 1$, the elements of $\Algebra^*$ and 
define a continuous and injective map $J:(X,\tau)\to\R^\infty$ by
\begin{equation*}
J(x):=\bigl(f_1(x), f_2(x), f_3(x),\ldots\bigr).
\end{equation*}
We fix $\tau$-compact sets $K_n\subset K_{n+1}$ and
$\mm(X\setminus K_n)\to 0$. The set 
$$
J^*:=\bigcup_{n=1}^\infty J(K_n)\subset J(X)
$$
is $\sigma$-compact in $\R^\infty$.

Defining $\nu_t\in\Probabilities{\R^\infty}$ by $\nu_t:=J_\#\mu_t$, $\cc:(0,T)\times\R^\infty\to\R^\infty$ by
$$
 \cc^i_t:=\begin{cases}
(\bb_t(f_i))\circ J^{-1}&\text{on $J^*$;}\\ \\ 0 &\text{otherwise,}
\end{cases}
$$
and noticing that
\begin{equation}\label{eq:boundcci}
|\cc^i_t|\circ J\leq |\bb_t|,\quad\text{$\mm$-a.e.\ in $X$,}
\end{equation}
the chain rule 
$$
\bb_t(\phi)(x)=\sum_{i=1}^n\frac{\partial\psi}{\partial z_i}( f_1(x),\ldots, f_n(x))\cc^i_t(x)
$$
for $\phi(x)=\psi(f_1(x),\ldots,f_n(x))$ 
shows that the assumption of Theorem~\ref{thm:superpoRinfty} are satisfied by $\nu_t$ with velocity $\cc$, because
\eqref{eq:boundcci} and $\mu_t\ll\mm$ give $|\cc^i_t|\leq |\bb_t|\circ J^{-1}$ $\nu_t$-a.e.\ in $\R^\infty$.
Notice also that all measures $\nu_t$ are concentrated on $J^*$.

As a consequence we can apply Theorem~\ref{thm:superpoRinfty} to obtain $\llambda\in\Probabilities{C([0,T];\R^\infty)}$
concentrated on solutions $\gamma\in AC_w([0,T];\R^\infty)$ to the ODE $\dot\gamma=\cc_t(\gamma)$ such that $(e_t)_\#\llambda=\nu_t$ for all $t\in [0,T]$. 
Since all measures $\nu_t$ are concentrated on $J^*$, one has
$$
\text{$\gamma(t)\in J^*$ for $\llambda$-a.e.\ $\gamma$, for all $t\in [0,T]\cap\Q$.}
$$
Denoting by $L$ the $\llambda$-negligible set above, the curve $\eta:=J^{-1}\circ\gamma:[0,T]\cap\Q\to X$ is defined for all
$\gamma\notin L$. For $s,\,t\in [0,T]\cap\Q$ with $s<t$, the estimate
\begin{eqnarray*}
\sup_i|f_i(\eta(s))-f_i(\eta(t))|&=&\sup_i|\gamma_i(s)-\gamma_i(t)|\leq \sup_i\int_s^r |\cc^i_r|(\gamma(r))\,dr\\&\leq& 
\int_s^r|\bb_r|(J^{-1}(\gamma(r))\chi_{J^*}(\gamma(r))\, dr
\end{eqnarray*}
shows that for $\llambda$-a.e.\ $\gamma$ one has $\gamma\in AC([0,T];\R^\infty)$ (by a density argument), 
and the corresponding curve $\eta$ is uniformly continuous in $[0,T]\cap\Q$ w.r.t.\ $d_{\Algebra^*}$.
By assumption (b) on $\Algebra^*$ it follows that $\eta$ has a unique extension to a continuous curve from $[0,T]$ to
$(X,\tau)$, and that $\gamma=J\circ\eta$ has image contained in $J(X)$. For these reasons, it makes sense to define 
$$
\eeta:=\Theta_\#\llambda
$$
where $\Theta:\{\gamma\in AC([0,T];\R^\infty): \gamma([0,T]\cap\Q)\subset J^*\}\to C([0,T];(X,\tau))$ is the $\llambda$-measurable
map $\gamma\mapsto \Theta(\gamma):=J^{-1}\circ\gamma$. Since
$(J^{-1})_\#\nu_t=\mu_t$, we obtain immediately that $(e_t)_\#\eeta=\mu_t$.

Let $i\geq 1$ be fixed. Since $f_i\circ\eta=p_i\circ\gamma$, taking the definition of $\cc_i$ into account we obtain that 
$f_i\circ\eta$ is absolutely continuous in $[0,T]$ and that
\begin{equation}\label{eq:partial_derivatives}
\text{$(f_i\circ\eta)'(t)=\bb_t(f_i)(\eta(t))$ for a.e.\  $t\in (0,T)$, for $\eeta$-a.e.\ $\eta$.}
\end{equation}
The proof is then completed by showing that \eqref{eq:partial_derivatives} extends from $\Algebra^*$ to all of $\Algebra$,
with the weaker requirement of $W^{1,1}(0,T)$ regularity of $f\circ\eta$.
\end{proof}

\subsection{Existence and uniqueness of regular Lagrangian flows}

In this section, we consider a Borel family of derivations $\bb=(\bb_t)_{t\in (0,T)}$ satisfying
\begin{equation}\label{eq:basicbb}
|\bb|\in L^1(L^2).
\end{equation}
Under the assumption that the continuity equation has uniqueness of solutions in the class
\begin{equation*}
{\mathcal L}_+:=\bigl\{u\in L^\infty(L^1_+\cap L^\infty_+):\ \text{$t\mapsto u_t$ is weakly continuous in $[0,T]$, in duality with $\Algebra$}\bigr\}
\end{equation*}
for any nonnegative initial datum $\bar u\in L^1\cap L^\infty(X,\mm)$, and existence of solutions in the class 
\begin{equation}\label{eq:morenarrow}
\bigl\{u\in{\mathcal L}_+:\ \|u_t\|_{\infty}\leq  C(\bb)\|\bar u\|_\infty\,\,\,\forall t\in [0,T]\bigr\},
\end{equation}
for any nonnegative initial datum $\bar u\in L^1\cap L^\infty(X,\mm)$,
we prove existence and uniqueness of the regular Lagrangian flow  $\XX$ associated to $\bb$. 
Here uniqueness is understood in the
pathwise sense, namely $\XX(\cdot,x)=\YY(\cdot,x)$ in $[0,T]$ for $\mm$-a.e.\ $x\in X$, whenever $\XX$ and $\YY$ are
regular Lagrangian flows relative to $\bb$.

Notice that the need for a class of functions as large as possible where uniqueness holds is hidden in the proof of 
Theorem~\ref{thm:nosplitting}, where solutions are built by taking the time marginals of suitable probability measures on curves and uniqueness leads to a non-branching
result.

This, together with Theorem~\ref{thm:well_PDE}, is the main result of these notes. Notice that
the well posedness of the continuity equation is ensured precisely by the assumptions on $\bb$ and
on the semigroup stated in Theorem~\ref{thm:well_PDE}. However we preferred to state the
theorem in an abstract form, assuming a priori the validity of this property.

\begin{theorem}[Existence and uniqueness of the regular Lagrangian flow]\label{thm:uniflow}
Assume \eqref{eq:basicbb}, assumptions (a), (b) on $\Algebra^*$, 
and that the continuity equation induced by $\bb$ has uniqueness of solutions in ${\mathcal L}_+$ for all
nonnegative initial datum $\bar u\in L^1\cap\Lbm\infty$, as well as existence of solutions in the class \eqref{eq:morenarrow} for all
nonnegative initial datum $\bar u\in L^1\cap\Lbm\infty$. Then there exists a unique regular Lagrangian flow relative to $\bb$.
\end{theorem}
\begin{proof} (Existence) Let us build first a ``generalized'' flow. To this aim, we take $\bar u\equiv 1$ as initial datum and we apply first the assumption on existence of a solution
$u\in{\mathcal L}_+$ starting from $\bar u$, with $u_t\leq C(\bb)$, and then
the superposition principle stated in
Theorem~\ref{thm:superpo} to obtain $\eeta\in\Probabilities{C([0,T];(X,d_{\Algebra^*}))}$ whose time marginals are $u_t\mm$, concentrated
on solutions to the ODE $\dot\eta=\bb_t(\eta)$. Then, Theorem~\ref{thm:nosplitting} below (which uses the uniqueness
part of our assumptions relative to the continuity equation) provides a representation
$$
\eeta= \int_X \delta_{\eta_x} d\mm(x),
$$
with $\eta_x\in C([0,T];(X,d_{\Algebra^*}))$, such that $\eta_x(0)=x$ and $\dot\eta_x=\bb_t(\eta_x)$. 
Setting $\XX(\cdot,x)=\eta_x(\cdot)$, it follows that $\XX:X\times [0,T]$ is a regular flow, relative to $\bb$, 
since 
\begin{equation*}
\XX(t,\cdot)_\#(\bar u\mm)=(e_t)_\#\eeta=u_t\mm\leq C(\bb)\mm.
\end{equation*}

\noindent 
(Uniqueness) Given RLF's $\XX$ and $\YY$, consider the measure
$$
\eeta:=\frac{1}{2}\int_X \bigl(\delta_{\sXX(\cdot,x)}+\delta_{\sYY(\cdot,x)}\bigr)\,d\mm(x)
\in\Probabilities{C([0,T];(X,d_{\Algebra^*}))}.
$$
By applying Theorem~\ref{thm:nosplitting} below to $\eeta$, we obtain $\XX(\cdot,x)=\YY(\cdot,x)$ for $\mm$-a.e.\ $x\in X$.
\end{proof}

\begin{theorem}[No splitting criterion]\label{thm:nosplitting}
Assume \eqref{eq:basicbb} and that the continuity equation induced by $\bb$ has at most one solution in ${\mathcal L}_+$ for all
$\bar u\in L^1\cap L^\infty(X,\mm)$. Let $\eeta\in\Probabilities{C([0,T];(X,\tau))}$ satisfy:
\begin{itemize}
\item[(i)] $\eeta$ is concentrated on solutions $\eta$ to the ODE $\dot\eta=\bb_t(\eta)$;
\item[(ii)] there exists $L_0\in [0,\infty)$ satisfying
\begin{equation*}
(e_t)_\#\eeta\leq L_0\mm\qquad\forall \,t\in [0,T].
\end{equation*}
\end{itemize}
Then the conditional measures $\eeta_x\in\Probabilities{C([0,T];(X,\tau)})$ induced by the map $e_0$ are Dirac masses for
$(e_0)_\#\eeta$-a.e.\ $x$; equivalently, there exist $\eta_x\in C([0,T];(X,\tau))$ such that $\eta_x(0)=x$ and solving the ODE
$\dot\eta_x=\bb_t(\eta_x)$, satisfying
$\eeta= \int \delta_{\eta_x} d(e_0)_\#\eeta(x)$.
\end{theorem}

In order to prove that the conditional measures in Theorem~\ref{thm:nosplitting} are Dirac masses we use the following
simple criterion, whose proof can be achieved by using finer and finer partitions of the state space.

\begin{lemma}\label{dyadic}
Let $\eeta_x$ be a $\mm$-measurable family of positive finite measures in $C([0,T];(X,\tau))$ with the
following property: for any $t\in [0,T]$ and any pair of disjoint Borel sets $E,\,E'\subset X$
we have
\begin{equation}\label{martingale}
\eeta_x\left(\{\gamma:\ \eta(t)\in E\}\right)
\eeta_x\left(\{\eta:\ \eta(t)\in E'\}\right)=0
\quad\text{$\mm$-a.e.\ in $X$.}
\end{equation}
Then $\eeta_x$ is a Dirac mass for $\mm$-a.e.\ $x\in X$.
\end{lemma}

\begin{proof}
For a fixed $t\in (0,T]$ it suffices to check
that the measures $\lambda_x:=\eta(t)_\#\eeta_x$ are Dirac masses for $\mm$-a.e.\ 
$x$.  Then \eqref{martingale} gives $\lambda_x(E)\lambda_x(E')=0$ $\mm$-a.e.\ for
any pair of disjoint Borel sets $E,\,E'\subset X$. Let $\delta>0$ and let us
consider a partition of $X$ in countably many Borel sets $R_i$ having a diameter
less then $\delta$. Then, as $\lambda_x(R_i)\lambda_x(R_j)=0$ $\mm$-a.e.\ whenever
$i\neq j$, we have a corresponding decomposition of $\mm$-almost all of
$X$ in Borel sets $A_i$ such that ${\rm supp\,}\lambda_x\subset\overline{R}_i$
for any $x\in A_i$ (just take $\{\lambda_x(R_i)>0\}$ and subtract from it all
other sets $\{\lambda_x(R_j)>0\}$, $j\neq i$). 
Since $\delta$ is arbitrary the statement is proved.
\end{proof}

\noindent 
{\bf Proof of Theorem~\ref{thm:nosplitting}.} The heuristic idea is that if $\eeta_x$ are not Dirac masses, then a suitable localization
procedure provides, after taking the marginals, distinct solutions to the continuity equation with the same initial datum. In the proof
we use of course the transfer mechanism from the Lagrangian to the Eulerian side provided by Lemma~\ref{eq:fromODEtoPDE}.

If the thesis is false then $\eeta_x$ is not a Dirac mass
in a set of $\bar\mu$ positive measure and
we can find, thanks to Lemma~\ref{dyadic}, $t\in (0,T]$, disjoint Borel
sets $E,\,E'\subset X$ and a Borel set $C$ with $\mm(C)>0$ such that 
$$
\eeta_x\big(\{\eta:\ \eta(t)\in E\}\big)\eeta_x\big(\{\eta:\ \eta(t)\in E'\}\big)>0
\qquad\forall x\in C.
$$

Possibly passing to a smaller set having still strictly positive $\mm$-measure 
we can assume that 
\begin{equation}\label{useless}
0<\eeta_x(\{\eta:\ \eta(t)\in E\})\leq M\eeta_x(\{\eta:\ \eta(t)\in E'\})
\qquad\forall x\in C
\end{equation}
for some constant $M$. We define measures $\eeta^1,\,\eeta^2$ whose disintegrations
$\eeta^1_x$, $\eeta^2_x$ are given by (here and in the sequel we use $\nu\res B$ as an alternative notation for $\chi_B\nu$)
$$
\eeta_x^1:=\chi_C(x)\eeta_x\res\{\eta:\ \eta(t)\in E\},
\qquad
\eeta_x^2:=M\chi_C(x)\eeta_x\res\{\eta:\ \eta(t)\in E'\}
$$ 
and denote by $\mu_s^i$, $s\in [0,t]$, the solutions of the continuity equation
induced by $\eeta^i$.
Then
$$
\mu^1_0=\eeta_x(\{\eta:\ \eta(t)\in E\})\mm\res C,\qquad
\mu^2_0=M\eeta_x(\{\eta:\ \eta(t)\in E'\})\mm\res C,
$$
so that \eqref{useless} yields $\mu^1_0\leq\mu^2_0$. On the other hand,
$\mu_t^1$ is orthogonal to $\mu^2_t$: precisely, denoting by $\eeta_{tx}$ 
the image of $\eeta_x$ under the map $\eta\mapsto\eta(t)$, we have
$$
\mu_t^1=\int_C\eeta_{tx}\res E\,d\mu(x)\perp
M\int_C\eeta_{tx}\res E'\,d\mu(x)=\mu_t^2.
$$
In order to conclude, let $\rho:X\to [0,1]$ be the density of 
$\mu_0^1$ with respect to $\mu_0^2$ and set 
$$
\tilde{\eeta}_x^2:=M\rho(x)\chi_C(x)\eeta_x\res\{\gamma:\ \gamma(t)\in E'\}.
$$
We define the measure $\tilde{\eeta}^2$ whose disintegration
is given by $\tilde{\eeta}_x^2$ and denote by $\tilde{\mu}_s^2$, $s\in [0,t]$, 
the solution of the continuity equation induced by $\tilde{\eeta}^2$.

Notice also that $\mu^i_s\leq\mu_s$ and so $\mu_s^i\in {\mathcal L}_+$,
and since $\tilde{\eeta}^2\leq\eeta^2$ we obtain that $\tilde{\mu}_s^2\in {\mathcal L}_+$
as well. By construction $\mu_0^1=\tilde{\mu}^2_0$, while $\mu_t^1$ is orthogonal
to $\mu_t^2$, a measure larger than $\tilde{\mu}^2_t$. We have thus built two
different solutions of the continuity equation with the same initial condition, a contradiction.

\section{From transport to diffusion operators}

In this section, we describe some recent results contained in the second author's PhD thesis \cite{Trevisan-14}, where existence and uniqueness problems are investigated for diffusion operators, in the same abstract framework as above; two papers \cite{Trevisan-15a} and \cite{Trevisan-15b}, the first one dealing with
refined results in Euclidean spaces and the second one dealing with the abstract setting are in preparation. Motivations for such an extension come at least from two sides: first, the theory of infinitesimally Hilbertian metric measure spaces, mentioned in Example \ref{ex:infinitesimally-hilbertian} and currently under intensive development, requires new calculus tools, and diffusion processes, strongly connected with parabolic partial differential equations, provide a natural extension of flows and transport equations techniques. From another side, by studying diffusion processes, we are in a position to compare our results with other approaches in stochastic analysis, such as the theory of non-symmetric Dirichlet forms \cite{Ma-Rockner-92, Stannat-99} or that of singular diffusions developed e.g.\ in \cite{Eberle-99}: as one might expect, since it emerges as a generalization of the deterministic case, the DiPerna-Lions approach allows for dealing with possibly degenerate diffusions, which are not fully covered by these techniques.

Before we address the theory in the abstract setting, we recall that solutions of stochastic differential equations (SDE's)
\begin{equation}\label{eq:intro-sde} d x_t = \bb(t,x_t) dt + \ssigma(t,x_t) dW_t, \quad t \in (0,T), \, x_t \in \R^d, \end{equation}
can be understood in at least two different ways (here, besides the vector field $\bb: (0,T) \times \R^d \to \R^d$, we let $\ssigma: (0,T) \times \R^d \to \R^{d\times d}$ and $(W_t)_{t \in [0,T]}$ be a $d$-dimensional Wiener process, or Brownian motion). A first notion is that of ``strong'' solution, originally due to K.\ It\^o, where  $(x_t)_{t \in [0,T]}$ is an adapted functional of the Wiener process (and involves stochastic integration in  the It\^o sense). A second, weaker,  notion is related to the  martingale problem approach, classically developed by Stroock and Varadhan, see \cite{Stroock-Varadhan-06}. The initial observation is that, given any strong solution to \eqref{eq:intro-sde} and an ``observable'' $f \in C^2_c(\R^d)$, an application of It\^o formula entails that
\begin{equation} \label{eq:strong-solution} [0,T] \ni t\quad \mapsto \quad M^f_t := f(x_t) -f(x_0) - \int_0^t (\cL_s f) (x_s)ds \bra{ = \int_0^t \nabla f(x_s)  
\ssigma(s, x_s) dW_s } \end{equation}
is a martingale with respect to the Brownian filtration, and a fortiori with respect to the natural filtration of $(x_t)_{t \in [0,T]}$, i.e., the conditional expectation of $M^f_t$ given $(x_r)_{r\le s}$ is $M^f_s$, for every $s$, $t \in [0,T]$, with $s \le t$. In the expression above, we let $\cL$ be the Kolmogorov operator
\begin{equation}\label{eq:intro-diffusion-operator}  \cL_t f (x) := \sum_{i=1}^d \bb^i(t, x) \frac{\partial f}{\partial x^i}(x) + \frac{1}{2} \sum_{i,j=1}^d  \aa^{i,j}(t,x) \frac{\partial^2f}{\partial x^i \partial x^j}(x), \quad \text{for $(t,x) \in (0,T)\times \R^d$,} \end{equation}
and $\aa^{i,j}: = \sum_{k=1}^d \ssigma^{i,k} \ssigma^{j,k}$, which is regarded as an infinitesimal covariance.  By definition, the martingale problem consists in finding some process $(x_t)_{t \in [0,T]}$ such that $(M^f_t)_{t \in [0,T]}$ defined in \eqref{eq:strong-solution} is a martingale for every $f \in C^2_c(\R^d)$; moreover, since we  are focusing on processes with continuous paths, it turns out that one can always restate the martingale problem in terms of the law of $(x_t)_{t \in [0,T]}$ in the space $C([0,T];\R^d)$, i.e., a solution of the martingale problem is equivalently described by a probability measure $\eeta$ on $C([0,T];\R^d)$, and by letting the ``evaluation'' process $e_t: \eta \mapsto \eta(t)$ play the role of $x_t$.  Together with the introduction of this notion, Stroock and Varadhan were able to show well-posedness for martingale problems whenever the Kolmogorov operator is the sum of a bounded measurable vector field $\bb$ and a uniformly bounded, continuous and elliptic $\aa$ \cite[Theorem 7.2.1]{Stroock-Varadhan-06}. Rigorous links between the strong notion and martingale problem are provided by the classical Yamada-Watanabe theorem \cite{Yamada-Watanabe-71}.

Since the seminal works by Stroock and Varadhan, the theory of martingale problems has been developed, showing strong connections with Markov semigroups and PDE's, also in abstract frameworks, see e.g.\ \cite{Ethier-Kurtz-86}. Figalli \cite{Figalli-08} was the first to draw a precise connection between DiPerna-Lions theory and that of martingale problems, obtaining in particular well-posedness for a wide class of multi-dimensional diffusion processes, whose Kolmogorov operators $\cL$ are not necessarily continuous or elliptic, provided that some Sobolev regularity holds. More precisely, the main object of his study are the so-called Stochastic Lagrangian Flows, i.e.,\ Borel families $(\eeta(x))_{x \in \R^d}$ of probabilities measures on $C([0,T];\R^d)$, such that
\begin{enumerate}[(a)]
\item for $\scrL^d$-a.e.\ $x \in \R^d$, $\eeta(x)$ solves the martingale problem associated to $\cL$, starting from $x$;
\item the push-forward measures $(e_t)_{\sharp} \int \eeta(x)  d\scrL^d(x)$ (recall the definition of $e_t$ above) are absolutely continuous with respect to $\scrL^d$, with uniformly bounded densities.
\end{enumerate}
Let us stress the fact that, as in the deterministic theory, uniqueness is understood for flows, thus in a selection sense, and not simply for $\scrL^d$-a.e.\ initial datum.

Although the two conditions above are in formal correspondence with the deterministic case, Stochastic Lagrangian Flows are not necessarily (neither expected to be) deterministic functions of the initial datum only, the occurring uncertainty being ``encoded'' in the fact that $\eeta(x)$ is a probability measure. Moreover, a discrepancy with the deterministic theory occurs when the Kolmogorov operator $\cL$ reduces to a derivation, i.e.\ when $\aa=0$: a solution to the martingale problem may be non-trivially concentrated on possibly non-unique solutions of the ODE, so we may not recover regular Lagrangian flows. Despite this gap, Figalli's formalism of Stochastic Lagrangian Flows already proved to be an efficient tool to deal with SDE's under low regularity assumptions in Euclidean spaces and, together with the more PDE oriented paper \cite{Lebris-Lions-08}, has become the starting point for further developments.

\subsection{Diffusion operators}

To study diffusion processes in the measure space setting, the first notion that we introduce is the analogue of a diffusion operator \eqref{eq:intro-diffusion-operator}, as an operator $\Algebra \ni f \mapsto \cL (f) \in \Lbm 0$. Here, the main difficulty is that second order regularity seems to be required to formulate an analogue of  \eqref{eq:intro-diffusion-operator}. Although a theory of second order calculus has been recently introduced in metric measure spaces whose Riemannian Ricci curvature is bounded from below \cite{gigli-WIP}, for many purposes it will suffice a formulation modelled on $\Gamma$-calculus, based on the validity of Leibniz rule with an extra term,
\[ \frac 1 2 \sqa{ \cL(f g ) - f \cL(g) - g \cL(f)} = \aa(f, g), \quad \text{for every $f,\,g \in \Algebra$.}\]
The (first order) bilinear operator appearing in the right hand side above can be generalized as an abstract bilinear map $(f,g) \mapsto \aa(f,g)$ on $\Algebra^2$. More precisely, one requires that for fixed $g$, the map $f \mapsto \aa(f, g)$ is a derivation in the sense of Definition~\ref{def:derivation} and that, for every 
$f,\,g \in \Algebra$, $\aa(f,g) = \aa(g,f)$, as well as $\lambda \Gamma(f) \le \aa(f,f) \le \Lambda \Gamma(f)$, for some non-negative functions $\lambda$, $\Lambda \in L^0$. If $\lambda$ is uniformly greater than some constant $\lambda_0 >0$ we say that $\cL$ is elliptic; a somewhat opposite case is that of a derivation $\cL(f) = \bb(f)$, where $\aa =0$.

The approach to diffusion operators sketched above can be made rigorous and already provides some non-trivial existence results, as well as abstract Eulerian/Lagrangian correspondence (via a suitable superposition principle), with the main drawback that formulation of density assumptions on $\cA$ are cumbersome, as well as uniqueness results, relying on regularity assumptions on the coefficients. For the sake of simplicity, in what follows  we study an explicit example, namely a diffusion operator in the form
\begin{equation}\label{eq:model-diffusion} \cL_t f := \bb_t(f) + \frac 1 2 a_t \Delta f, \quad \text{for $f \in \Algebra$,}\end{equation}
where $\bb$ is a Borel time-dependent family of derivations, $a: (0,T) \times X \to [0,\infty)$ is Borel, and  we assume that $\Algebra \subseteq \V_\infty$ is also included in the space of functions $f$ such that $\Delta f \in L^2\cap L^\infty(X,\mm)$. The chain rules for derivations and the Laplacian entail
\[ \frac 1 2  \sqa{ \cL (f   g ) - g \cL(f) - f\cL(g)} = a \Gamma(f,g), \quad \text{for every $f$, $g \in \Algebra$,}\]
showing that $a$ should be non-negative.

Regularity assumptions on the ``coefficients'' of $\cL$ can be then easily stated, in terms of assumptions on $\bb$, $a$ and the gradient vector field $\nabla a: f \mapsto \Gamma(f, a)$. Below, to show uniqueness in the (possibly) non-elliptic case, we impose bounds on the symmetric part of the latter vector field, i.e., on the already quoted Bakry's Hessian of the function $a$ \cite{Bakry-94}.

\subsection{Fokker-Planck equations}

In this section, we sketch how the Eulerian side of the theory  generalizes to diffusion operators in the form \eqref{eq:model-diffusion}.  We study the analogue of the continuity equation, given by the Fokker-Planck equation
\[
\frac{d}{dt}u_t = \cL^*_t u_t, \quad \text{ on $(0,T) \times X$,}
\]
defined in duality with $\Algebra$, under the integrability assumptions
$$
u\in L^\infty(L^\infty), \qquad |\bb|\in L^1(L^{1}),\qquad a \in L^1(L^{1}).
$$
We always consider the following weak formulation:
\begin{equation}\label{eq:very-weak-fpe}
\frac{d}{dt}\int_X u_t f\,d\mm=\int_X \cL_t f d\mm = \int_X\sqa{  \bb_t( f) + \frac 1 2 a_t \Delta f} u_td\mm\quad\text{in $\mathcal D'(0,T)$, for all $f\in\Algebra$.}
\end{equation}
As in the case of continuity equation, a density argument shows that there exists a unique representative $t\mapsto u_t$ which is continuous in the duality with $\Algebra$, and we may also encode the initial condition $u_0 = \bar{u}$, and the mass conservation property conservation follows from the fact that $\cL 1 = 0$. The following existence result can be proved arguing as in Theorem \ref{thm:existence}, the proof being easier because functions in $\Algebra$ are more regular than the deterministic case, in particular we assume $\Delta f \in L^\infty$ for $f \in \Algebra$.

\begin{theorem} Let $\bar u\in L^\infty(X,\mm)$. Assume that $|\bb|$, $a\in L^1(L^2)$,
$\div\bb - \Delta a/2\in L^1(L^2)$ and that its negative part $(\div\bb - \Delta a/2)^-$ belongs to $L^1(L^\infty)$. Then, there
exists a solution $u\in L^\infty(L^\infty)$ to \eqref{eq:CEbis}. In addition, we can build $u$ in such a way that
$u\geq 0$ if $\bar u\geq 0$.
\end{theorem}

In case $\cL$ is elliptic, one can also obtain a uniqueness result for solutions belonging to the space $L^2((0,T);\V)$. Indeed, for $u$ belonging to this class and assuming e.g.\ that $a$, $|\bb|$, $|\nabla a| \in L^\infty(L^\infty)$, we may integrate by parts in \eqref{eq:very-weak-fpe}, obtaining 
\[\frac{d}{dt}\int_X u_t f\,d\mm = \int_X\sqa{  \bb_t( f) - \frac 1 2\Gamma( a_t , f)}u_t d\mm - \frac 1 2 \int_X a_t \Gamma(u_t, f) d\mm\quad\text{in $\mathcal D'(0,T)$, for all $f\in\Algebra$.}\]
which is actually a more standard weak formulation of the PDE (while \eqref{eq:very-weak-fpe} is sometimes referred in the literature as a very weak formulation). It follows that $\partial_t u \in L^2(\V^*)$ so we may let $u$ as a test function, and obtain the usual energy estimate
\[\frac{d}{dt}\int_X u_t^2 d\mm= \int_X\sqa{  2\bb_t(u_t) - \Gamma( a_t , u_t)}u_t d\mm - \int_X a_t \Gamma(u_t, u_t) d\mm,\]
which, by ellipticity and straightforward bounds such as $2\abs{\bb_t(u_t)} \le \veps^{-1} |\bb_t|^2 u_t^2 +  \veps\Gamma(u_t)$ for $\veps$ small enough, entails the differential inequality $\frac{d}{dt} \nor{ u_t}_2^2 \le c(t)  \nor{u_t}_2^2$, for some $c \in L^1(0,T)$. Eventually uniqueness follows by taking the difference between two solutions and applying Gronwall lemma. Let us also notice that, if we look for uniqueness in $L^\infty(L^\infty) \cap L^2(\V)$, the assumption $|\nabla a| \in  L^\infty$ can be relaxed to $|\nabla a| \in L^2(L^2)$, if we integrate by parts the term $\Gamma( a_t , u_t) u_t$ and impose a uniform bound on $\Delta a$.

In any case, uniqueness results for solutions belonging to spaces where we impose regularity bounds w.r.t.\ the spatial variables, e.g.\ $\Gamma(u) \in L^1(L^1)$ or $a\Gamma(u) \in L^1(L^1)$ (studied e.g.\ in \cite{Lebris-Lions-08} in Euclidean spaces, under low regularity assumptions on the coefficients) are not well suited for translation to the Lagrangian side of the theory, the main issue being the  proof of Theorem~\ref{thm:nosplitting}, which does not seem to be
compatible with such regularity requirements. Therefore, we study uniqueness in spaces such as $L^\infty(L^\infty)$; we give two well-posedness statements, one for possibly degenerate diffusions and the other one for the elliptic case.

\begin{theorem}\label{thm:wp-fpe-degenerate} Assume that $\Algebra$ is dense in $\V_4 \cap D(\Delta)$ and that $L^4$-$\Gamma$ inequalities hold. If
$$
|\bb|,\, |D^{sym}\bb_t|\in L^1(L^2),\qquad
\div\bb_t\in L^1(L^\infty)
$$
and 
$$
a, \sqrt{\Gamma(a)}, |\operatorname{Hess}[a]|\in L^1(L^2),\qquad
\Delta a\in L^1(L^\infty),
$$
then the Fokker-Planck equation \eqref{eq:CE} with initial datum $\bar u\in L^\infty(X,\mm)$ is well posed in $L^\infty(L^\infty)$.
\end{theorem}

Roughly speaking, we assume above first order Sobolev regularity for the drift, and second order regularity for the diffusion coefficient. In case of bounded elliptic diffusions, well-posedness holds even if we ``remove'' one order of regularity; however, we impose Lipschitz regularity with respect to $t$ for the map $t \mapsto a(t,x)$.

\begin{theorem} \label{thm:wp-fpe-elliptic}Assume that $\Algebra$ is dense in $\V_2 \cap D(\Delta)$ and that $a \ge \lambda$, a.e.\ on $(0,T)\times X$, for some $\lambda >0$. If
$$
|\bb|,\, a,\, \partial_t a \in L^\infty(L^\infty), \qquad
\sqrt{\Gamma(a)} \in L^1(L^2),\qquad
\Delta a\in L^1(L^\infty),
$$
then the Fokker-Planck equation \eqref{eq:CE} with initial datum $\bar u\in L^\infty(X,\mm)$ is well posed in $L^\infty(L^\infty)$.
\end{theorem}

Both theorems are proved by means of a smoothing scheme akin to that of Theorem \ref{thm:well_PDE}: in the former, we regularize using the semigroup $P$, so the commutator between $P$ and $\cL$ appears; in the latter, we regularize  by means of the semigroup associated to the Dirichlet form $L^2(L^2) \ni f \mapsto \int_0^T \int_X a_t \Gamma(f_t)d\mm dt$, so its commutator with $\partial_t$ appears (this is essentially the reason why we require $t \mapsto a_t$ to be Lipschitz). Commutators are then estimated using a similar interpolation and duality argument as in Theorem \ref{thm:commutator}, but in the case of possibly degenerate diffusions, we estimate the commutator between $P$ and $\cL$ relying on a somehow more precise
``second order'' interpolation, which eventually leads to terms involving the Hessian of $a$.

\subsection{Martingale problems}

As in the deterministic case, the Lagrangian side of the theory relies on a suitable superposition principle, lifting any weak solution of the Fokker-Planck equation to some solution of the correspondent martingale problem, and a ``no splitting criterion'' (akin to Theorem \ref{thm:nosplitting}), transferring uniqueness for PDE's to uniqueness for processes. Before we sketch the stochastic counterparts of these results, we rigorously define solutions of the martingale problem: as in Section~\ref{sec:ODE-RLF}, we fix a countable set $\Algebra^* \subseteq \Algebra$ (consisting of $\tau$-continuous functions), but we limit ourselves to processes with continuous paths with respect to coarsest topology $\tau_0$ on $X$ which makes continuous every $f \in \Algebra^*$.

\begin{definition}[regular solution to the martingale problem]
Let $|\bb|$, $a \in L^1(L^1)$. A probability measure $\eeta \in \Probabilities{ C([0,T]; (X,\tau_0) )}$ is a regular solution to the martingale problem 
associated to $\cL( = \bb + \frac 1 2 a \Delta)$ if one has $(e_t)_\sharp \eeta \in L^\infty(L^\infty)$ and, for every $f \in \Algebra^*$, the process
\[ [0,T] \ni t \mapsto M^f_t := f( e_t ) - f( e_0 )- \int_0^t (\cL_s f) ( e_s )ds\]
is a martingale with respect to the natural filtration associated to $(e_t)_{t \in [0,T]}$.
\end{definition}

The process $M^f$ is seen to be integrable, adapted and continuous, so that the martingale property reduces to show that the conditional expectation (w.r.t.\ $\eeta$) of $M^f_t$ given $(e_r)_{r\le s}$ coincides with $M^f_s$, for every $s$, $t \in [0,T]$, with $s \le t$. 

It is not difficult to show that the quadratic variation process of the martingale $M^f$ is $t \mapsto \int_0^t a_s( e_s) \Gamma(f)(e_s) ds$ (which would be a trivial consequence of the It\^o integral representation in \eqref{eq:strong-solution}). As a consequence, when $a=0$ and $\cL = \bb$ is a derivation, $M^f$ is a constant martingale and regular solutions of the martingale problem are concentrated on solutions to the ODE $\dot\eta=\bb_t(\eta)$, in the sense of Definition \ref{def:ourODE}. 
Another, less trivial,  consequence, essentially based on Burkholder-Gundy inequalities, is an estimate on the modulus of continuity of $t \mapsto M^f_t (\eta)$, for $\eeta$-a.e.\ $\eta \in C([0,T]; (X,\tau_0) )$ in terms of the integrability of $\abs{\bb}$ and $a$, which we may regard as an analogue of \eqref{eq:inequality-or-equality}, but of course $M^f$ is in general only H\"older continuous (for brevity, we omit to state the precise result). Let us also point out that the choice of $\tau_0$ in place of $\tau_{\Algebra^*}$, that of the deterministic case, is motivated by the apparent difficulty to provide estimates on the modulus of continuity of the process with respect to $\sup$-type distances on $X$, such as $d_{\Algebra^*}$, an application of It\^o formula would involve second-order derivatives of the function $(z_i)_{i\ge 1} \mapsto \sup_{i\ge 1} \abs{z_i}$.

With this notation, the superposition principle for diffusion processes in metric measure spaces reads as an analogue of Theorem~\ref{thm:superpo}, where we ``lift'' any narrowly continuous solution $(u_t \mm)_{t \in [0,T]} \subset\Probabilities{X}$ of the Fokker-Planck equation to some solution $\eeta$ of the correspondent martingale problem, with $1$-marginals $(e_t)_\sharp \eeta = u_t \mm$, for $t \in [0,T]$. Its proof is based on a similar reduction first to the case $X = \R^\infty$ and then to $X = \R^d$, although one has to fill several technical aspects.

We conclude by pointing out that, in order to transfer well-posedness from the Eulerian to the Lagrangian description, we can not use a direct analogue of 
Theorem~\ref{thm:nosplitting}: we actually expect some ``good'' splitting to occur (notice also that in the proof of we are conditioning at time $t = 0$, with respect to a ``future event'', namely $\cur{\eta: \eta(t) \in E}$, which would give rise to some problem also for strong solutions). To obtain a similar conclusion, we rely instead on a weaker criterion, originally due to Stroock and Varadhan \cite[Theorem 6.2.3]{Stroock-Varadhan-06} and then modified in the framework of finite-dimensional regular solutions in \cite[Proposition 5.5]{Figalli-08}: its proof in the abstract framework requires only minor changes.

\begin{lemma}
Assume that, for every $s \in [0,T]$ and every probability density $\overline u \in L^\infty(X,\mm)$, there exists a unique solution $(u_t\mm)$ 
in $C([s,T]; \Probabilities{X}) \cap L^\infty(L^\infty)$ of the Fokker-Planck equation
\[ \partial_t u_t = \cL_t^* u_t, \quad \text{in $(s,T)\times \R^d$, with $u_s = \overline{u}$}.\]
Then, for every $s \in [0,T]$ and every probability density $\overline u \in L^\infty(X,\mm)$, there exists a unique regular solution  $\eeta$ of the correspondent martingale problem on $C([s,T];(X, \tau_0))$ with $(e_s)_\sharp \eeta = \overline u \mm$.
\end{lemma}

Under the assumptions of Theorem~\ref{thm:wp-fpe-degenerate} or Theorem~\ref{thm:wp-fpe-elliptic}, we are then in a position to deduce well-posedness for martingale problems, since these hold also on intervals $[s,T]$, for $s \in [0,T]$. A disintegration argument with respect to the initial condition, as in Theorem~\ref{thm:uniflow}, eventually leads to well-posedness of suitable versions of Stochastic Lagrangian Flows in abstract measure spaces.

\def\cprime{$'$}
\providecommand{\bysame}{\leavevmode\hbox to3em{\hrulefill}\thinspace}
\providecommand{\MR}{\relax\ifhmode\unskip\space\fi MR }
\providecommand{\MRhref}[2]{%
  \href{http://www.ams.org/mathscinet-getitem?mr=#1}{#2}
}
\providecommand{\href}[2]{#2}

\end{document}